\documentclass{statsoc}
\usepackage{amssymb, amsmath, mathrsfs}
\usepackage[left=1in,right=1in,top=1in,bottom=1in]{geometry}

\newtheorem{assumption}{Condition}
\newtheorem{lemma}{Lemma}
\newtheorem{proposition}{Proposition}
\newtheorem{theorem}{Theorem}
\newtheorem{definition}{Definition}

\newcommand{\ba}{\mbox{\bf a}}
\newcommand{\bb}{\mbox{\bf b}}
\newcommand{\bd}{\mbox{\bf d}}

\newcommand{\bff}{\mbox{\bf f}}

\newcommand{\bu}{\mbox{\bf u}}
\newcommand{\bv}{\mbox{\bf v}}
\newcommand{\bw}{\mbox{\bf w}}
\newcommand{\bx}{\mbox{\bf x}}
\newcommand{\by}{\mbox{\bf y}}
\newcommand{\bz}{\mbox{\bf z}}
\newcommand{\bA}{\mbox{\bf A}}
\newcommand{\bB}{\mbox{\bf B}}
\newcommand{\bC}{\mbox{\bf C}}

\newcommand{\bG}{\mbox{\bf G}}

\newcommand{\bH}{\mbox{\bf H}}

\newcommand{\bV}{\mbox{\bf V}}

\newcommand{\bX}{\mbox{\bf X}}

\newcommand{\bY}{\mbox{\bf Y}}
\newcommand{\bZ}{\mbox{\bf Z}}
\newcommand{\bone}{\mbox{\bf 1}}
\newcommand{\bzero}{\mbox{\bf 0}}
\newcommand{\bveps}{\mbox{\boldmath $\varepsilon$}}

\newcommand{\bbeta}{\mbox{\boldmath $\beta$}}
\newcommand{\bdelta}{\mbox{\boldmath $\delta$}}
\newcommand{\btheta}{\mbox{\boldmath $\theta$}}
\newcommand{\bgamma}{\mbox{\boldmath $\gamma$}}
\newcommand{\bPsi}{\mbox{\boldmath $\Psi$}}

\newcommand{\bmu}{\mbox{\boldmath $\mu$}}

\newcommand{\hbbeta}{\widehat\bbeta}
\newcommand{\hbgamma}{\widehat\bgamma}

\newcommand{\var}{\mathrm{var}}
\newcommand{\cov}{\mathrm{cov}}

\newcommand{\Sig}{\mathbf{\Sigma}}

\newcommand{\tr}{\mathrm{tr}}
\newcommand{\diag}{\mathrm{diag}}

\newcommand{\supp}{\mathrm{supp}}

\def\t{^T}
\def\toD{\overset{\mathscr{D}}{\longrightarrow}}

\title[Model Selection Principles in Misspecified Models]
{Model selection principles in misspecified models
\thanks{Lv's research was supported by NSF CAREER Award DMS-0955316 and Grant DMS-0806030 and 2008 Zumberge Individual Award from USC's James H. Zumberge Faculty Research and Innovation Fund. Liu's research was supported by NSF Grant DMS-0706989 and NIH Grant R01-HG02518-02. We sincerely thank the Joint Editor, an Associate Editor, and a referee for their constructive comments that significantly improved the paper. The first draft of the paper was written in 2008.}
}

\author{Jinchi Lv}
\address{University of Southern California, Los Angeles, USA}
\email{jinchilv@marshall.usc.edu}
\author[J. Lv and J. S. Liu]{and Jun S. Liu}
\address{Harvard University, Cambridge, USA}
\email{jliu@stat.harvard.edu}

\begin{document}

\begin{abstract}
Model selection is of fundamental importance to high-dimensional modeling featured in many contemporary applications. Classical principles of model selection include the Bayesian principle and the Kullback-Leibler divergence principle, which lead to the Bayesian information criterion and Akaike information criterion, respectively, when models are correctly specified. Yet model misspecification is unavoidable in practice. We derive novel asymptotic expansions of the two well-known principles in misspecified generalized linear models, which give the generalized BIC (GBIC) and generalized AIC (GAIC). A specific form of prior probabilities motivated by the Kullback-Leibler divergence principle leads to the generalized BIC with prior probability ($\mbox{GBIC}_p$), which can be naturally decomposed as the sum of the negative maximum quasi-log-likelihood, a penalty on model dimensionality, and a penalty on model misspecification directly. Numerical studies demonstrate the advantage of the new methods for model selection in both correctly specified and misspecified models.

\end{abstract}

\section{Introduction} \label{sec1}
We consider in this article the   type of data $(y_i, x_{i1}, \cdots,
x_{ip})_{i = 1}^n$, where the $y_i$'s are
independent observations of the response variable $Y$ conditional
on its
covariates, or explanatory variables, $(x_{i1}, \cdots,
x_{ip})\t$. When the dimensionality or the number of covariates $p$ is
large compared with the sample size $n$, it is  desirable to
produce sparse models that involve only a small subset of predictors. With such models one can improve the prediction accuracy
and enhance the model interpretability. See, for example, Fan and Lv (2010) for an
overview of recent progresses in high-dimensional variable selection
problems. A natural and fundamental problem that arises from such a
task is how to
compare models with  different sets of predictors.

Two classical principles of model selection are the Kullback-Leibler (KL)
divergence principle and the Bayesian principle,
which lead to the Akaike information criterion (AIC) by Akaike (1973, 1974) and
Bayesian information criterion (BIC) by Schwartz (1978), respectively,
when the models are correctly specified. The AIC and BIC have been proved to be powerful tools for model selection in various settings; see, for example, Burnham
and Anderson (1998) for an account of these developments. Stone (1977)
showed heuristically the asymptotic equivalence
of AIC and the leave-one-out cross-validation. Bozdogan (1987) studied asymptotic properties of the
AIC and  BIC, showing that in parametric settings, AIC asymptotically has a positive
probability of overestimating the true dimension, while BIC is
asymptotically consistent in estimating the true model. The model selection consistency of BIC has also been shown in Casella
et al. (2009). Hall (1990) compared AIC with KL
cross-validation in the setting of
histogram density estimation. Yang and Barron (1998) studied the
asymptotic property of model selection criteria related to the AIC and
minimum description length principles in density estimation. Other
work on model selection includes the risk inflation criterion (Foster
and George, 1994), the generalized information criterion (Konishi and
Kitagawa, 1996), Bayesian measure using deviance information criterion
(Spiegelhalter et al., 2002; Gelman et al., 2004),
model evaluation using the absolute prediction error (Tian
et al., 2007), tuning parameter selection in
penalization method (Wang, Li and Tsai, 2007; Fan and Tang, 2012), the parametricness index
(Liu and Yang, 2011), and many extensions of the AIC and BIC (see,
for example, Bozdogan, 1987 and 2000; Bogdan, Ghosh and Doerge, 2004; Chen and Chen, 2008; \.{Z}ak-Szatkowska and Bogdan, 2011).

In parametric modeling, one chooses \emph{a priori}
a family of distributions and then finds a model in this family that
best fits the data, often based on the maximum likelihood
principle. Despite this common practice, it has been a common wisdom in statistics that ``all models are wrong,
but some are more useful than others."  Thus, in a broad sense all
models are misspecified. Neither AIC nor BIC, however, explicitly accounts for such model misspecification.

When the model family is  misspecified, evaluating whether a model
selection procedure is good or bad becomes more subtle because we can no
longer simply say that a selected model is ``correct'' or
``incorrect.'' We can mathematically define that a procedure is
consistent if it can asymptotically select the model in the family that is theoretically
closest to the true model, according to a certain discrepancy measure
such as the KL divergence. In the same token, we can define a
procedure to be better than others if on average it can select a model closer to the
true model than others, according to a certain discrepancy measure. In
practice, this type of discrepancy computation can, however, be
nontrivial. We content ourselves with the scenario of generalized linear
models (GLMs) (McCullagh and Nelder, 1989) with $p$
predictors. In this case, a procedure is better than another if on
average it can select more relevant predictors than the
other. Prediction is also another important measure for model evaluation.

In this paper, we derive  asymptotic expansions of the Bayesian
and KL divergence principles in misspecified GLMs. The general idea
appears also applicable to other model settings.
The technical results lead us to propose the generalized BIC (GBIC)
and generalized AIC (GAIC). A specific choice of prior probabilities
for the models motivated by the KL divergence principle leads us to
propose the generalized BIC with prior probability ($\mbox{GBIC}_p$),
which can be written as the sum of the negative maximum
quasi-log-likelihood, a penalty on model dimensionality, and a penalty
on model misspecification directly. Our results can also be viewed as an extension
of the classical approaches to non-independent and identically distributed (non-i.i.d.) settings.

The rest of the paper is organized as follows. In Section \ref{sec2},
we discuss the quasi-maximum likelihood estimation (QMLE) in misspecified
GLMs and two classical model selection principles. We investigate
asymptotic expansions of these principles in misspecified models and
present resulting model selection criteria in Section \ref{sec3}. In
Section \ref{sec4}, we provide all technical conditions and
investigate asymptotic properties of the quasi-maximum likelihood estimator (QMLE). We present several
numerical examples to illustrate the finite-sample performance of
these model selection criteria in both correctly specified and
misspecified models in Section \ref{sec5}. In Section \ref{sec6}, we
provide some discussions of our results and their implications as well
as some possible extensions. All technical details are relegated to the
Appendix.

\section{Model misspecification and model selection principles} \label{sec2}

\subsection{Quasi-maximum likelihood estimation} \label{sec2.1}
To facilitate the presentation, we introduce some notation and state
some basic facts of the QMLE in
misspecified GLMs with deterministic design matrices.
Conditional on the covariates, we assume that the
$n$-dimensional random response vector
$\bY = (Y_1, \cdots,
Y_n)\t$ has a true unknown distribution $G_n$ with density function
\begin{equation} \label{013}
g_n(\by) = \prod_{i = 1}^n g_{n, i}(y_i),
\end{equation}
where $\by = (y_1, \cdots, y_n)\t$. Each observation $y_i$ of the
response variable may depend on all or some of the $p$ predictors
$x_{ij}$. Model (\ref{013}) entails that all components of $\bY$
are independent but not necessarily identically distributed.

Consider an arbitrary subset $\mathfrak{M}
\subset \{1, \cdots, p\}$ with $d = |\mathfrak{M}| \leq n \wedge
p$. Define $\bx_i = (x_{ij}: j \in \mathfrak{M})\t$ and $\bX = (\bx_1,
\cdots, \bx_n)\t$, an $n \times d$ deterministic design matrix. Since
the true model $G_n$ is unknown, we
choose a family of generalized linear models $F_n(\cdot,\bbeta)$ as
our working models,  with density function
\begin{equation} \label{003}
f_n(\bz, \bbeta) d\mu_0(\bz)= \prod_{i = 1}^n f_0(z_i, \theta_i)
d\mu_0(z_i) \equiv \prod_{i = 1}^n \exp \left[\theta_i z_i -
  b(\theta_i)\right] d\mu(z_i),
\end{equation}
where $\bz = (z_1, \cdots, z_n)\t$, $ \btheta = (\theta_1, \cdots,
\theta_n)\t = \bX \bbeta $
with $\bbeta \in \mathbb{R}^d$, $b(\theta)$ is a smooth convex
function, $\mu_0$ is the Lebesgue measure on Euclidean spaces, and $\mu$ is some fixed measure on $\mathbb{R}$. Clearly $\{f_0(z, \theta): \theta \in \mathbb{R}\}$ is
a family of distributions in the regular exponential family and may
not contain $g_{n, i}$'s. Denote by
$\bb(\btheta) = (b(\theta_1), \cdots, b(\theta_n))\t$ and
$\bmu(\btheta) = (b'(\theta_1), \cdots,
b'(\theta_n))\t$, and define a  matrix-valued function $\Sig(\btheta) =
\diag\{b''(\theta_1), \cdots,
b''(\theta_n)\}$. It is a well-known fact that for any
$n$-dimensional random vector $\bZ$ with distribution $F_n(\cdot, \bbeta)$
given by (\ref{003}), we have
\begin{equation} \label{073}
E \bZ = \bmu(\bX \bbeta) \quad \text{and} \quad \cov(\bZ) =
\Sig(\bX \bbeta).
\end{equation}

Theorem \ref{thm4} in Section \ref{sec4.2} shows that under some regularity conditions, there exists a unique distribution
$F_n(\cdot, \bbeta_{n, 0})$
in the family of misspecified GLMs in (\ref{003}) that has the smallest KL
divergence from the true model $G_n$, where the KL divergence (Kullback
and Leibler, 1951) of the posited GLM
$F_n(\cdot, \bbeta)$ in (\ref{003}) from the true model $G_n$ is given by
\begin{equation} \label{008}
I(g_n; f_n(\cdot, \bbeta)) = \int [\log
  g_n(\bz)] g_n(\bz) d\bz - \int [\log f_n(\bz, \bbeta)] g_n(\bz)
d\bz  = \sum_{i =
  1}^n I(g_{n, i}; f_0(\cdot, \theta_i)),
\end{equation}
where $I(g_{n, i}; f_0(\cdot, \theta_i)) = \int [\log g_{n, i}(z)] g_{n, i}(z) dz - \int [\log f_0(z, \theta_i)] g_{n, i}(z) dz$. Thus, the parameter $\bbeta$ is identifiable in such misspecified models. This uniqueness entails that $F_n(\cdot,
\bbeta_{n, 0}) = G_n$ when the model is correctly
specified, that is, $G_n \in \{F_n(\cdot, \bbeta): \bbeta \in
\mathbb{R}^d\}$. The QMLE of the
$d$-dimensional parameter vector $\bbeta_{n, 0}$ is
defined as
\begin{equation} \label{005}
\hbbeta_n = \mbox{arg}\max_{\bbeta \in \mathbb{R}^d} \ell_n(\by, \bbeta),
\end{equation}
where
\begin{equation} \label{002}
\ell_n(\by, \bbeta) = \log f_n(\by, \bbeta) = \by\t \bX
  \bbeta - \bone\t \bb(\bX \bbeta) + \sum_{i = 1}^n \log
\frac{d\mu}{d\mu_0}(y_i)
\end{equation}
is the quasi-log-likelihood function. Clearly $\hbbeta_n$ is the
solution to the  score equation
\begin{equation} \label{007}
\bPsi_n(\bbeta) \equiv \frac{\partial \ell_n(\by, \bbeta)}{\partial
  \bbeta} = \bX\t [\by - \bmu(\bX \bbeta)] = \bzero,
\end{equation}
which becomes the normal equation $\bX\t \by = \bX\t \bX \bbeta$ in
the linear regression model.

Under some regularity conditions, Theorem \ref{thm5} in Section
\ref{sec4.2} shows that the QMLE $\hbbeta_n$ is a consistent estimator
of $\bbeta_{n, 0}$, and Theorem \ref{thm6} further establishes its
asymptotic normality. The asymptotic covariance matrix of $\hbbeta_n$
is $\bA_n^{-1} \bB_n \bA_n^{-1}$, where
\begin{equation} \label{089}
\bB_n\equiv \cov\left[\bPsi_n(\bbeta_{n, 0})\right] = \cov\left(\bX\t
\bY\right) = \bX\t \cov(\bY) \bX,
\end{equation}
in view of (\ref{011}), with $\cov(\bY)  = \diag\{\var(Y_1), \cdots, \var(Y_n)\}$ by the
independence assumption, and $\bA_n = \bA_n(\bbeta_{n,
  0})$ with
\begin{equation} \label{090}
\bA_n(\bbeta)\equiv \frac{\partial^2 I(g_n; f_n(\cdot,
  \bbeta))}{\partial \bbeta^2}
= -\frac{\partial^2 \ell_n(\by, \bbeta)}{\partial \bbeta^2}
= \bX\t \Sig(\bX \bbeta) \bX
\end{equation}
under the smoothness condition. Throughout the paper, the expectation and
covariance are taken with respect to the true
distribution $G_n$ unless specified otherwise.

In view of (\ref{073}),
$\bA_n(\bbeta)$ defined in (\ref{090}) is exactly the covariance matrix of $\bX\t \bY$ when
$\bY$ has distribution $F_n(\cdot, \bbeta)$, and thus $\bA_n$ is the
covariance matrix of $\bX\t \bY$ under the best misspecified GLM
$F_n(\cdot, \bbeta_{n, 0})$, whereas $\bB_n$ defined in (\ref{089}) is the covariance matrix
of $\bX\t \bY$ under the true model $G_n$. It is well known in classical
maximum likelihood theory that the identity $\bA_n = \bB_n$ holds when the model is
correctly specified, with  $\bA_n$ and $\bB_n$ being the Hessian and outer product forms of the Fisher information matrix, respectively. These two matrices play a pivotal role in
quasi-maximum likelihood estimation of misspecified models.

\subsection{Bayesian and Kullback-Leibler divergence principles} \label{sec2.2}
Suppose we have a set of competing models: $\left\{\mathfrak{M}_m: m = 1,
\cdots, M\right\}$. A  typical Bayesian model selection procedure is to first give nonzero
prior probability $\alpha_{\mathfrak{M}_m}$ on each  model
$\mathfrak{M}_m$, and then prescribe a prior distribution
$\mu_{\mathfrak{M}_m}$ for the parameter vector in the corresponding model. Assume
that the density function of $\mu_{\mathfrak{M}_m}$  is bounded in
$\mathbb{R}^{\mathfrak{M}_m}
 = \mathbb{R}^{d_m}$  and locally bounded away from
 zero throughout the domain. The Bayesian principle of model selection  is
to choose the most probable model \emph{a posteriori}. That is, to choose model  $\mathfrak{M}_{m_0}$ such that
$
m_0 = \arg \max_{m \in \{1, \cdots, M\}} S(\by, \mathfrak{M}_m; F_n),
$
where
\begin{equation}\label{083}
 S(\by, \mathfrak{M}_m; F_n) = \log \int \alpha_{\mathfrak{M}_m}
\exp\left[\ell_n(\by, \bbeta)\right] d\mu_{\mathfrak{M}_m}(\bbeta)
\end{equation}
with the log-likelihood $\ell_n(\by, \bbeta)$ as in
(\ref{002}) and the integral over $\mathbb{R}^{d_m}$. Fix an arbitrary
model $\mathfrak{M} = \mathfrak{M}_m$ and let $\hbbeta_n
= \hbbeta_{n, m}$ be the QMLE. Using the Taylor expansion,
Schwartz (1978) showed that
$S(\by, \mathfrak{M}; F_n)$ can be asymptotically expanded as $\ell_n(\by, \hbbeta_n) - \frac{\log n}{2} |\mathfrak{M}| +
\log \alpha_{\mathfrak{M}}$ with $\hbbeta_n$ the maximum likelihood estimator (MLE) of
$\bbeta$,  for the case of
i.i.d. observations with correctly specified regular
exponential family of models. This expression allowed Schwartz to introduce the
seminal BIC for model selection:
\begin{equation} \label{063}
\mathrm{BIC}(\by, \mathfrak{M}) = -2 \ell_n(\by, \hbbeta_n) +
\left(\log n\right)
|\mathfrak{M}|,
\end{equation}
where a factor of 2 is included to make it consistent with the usual definition of AIC give in (\ref{062}). Schwartz's original arguments were generalized later by
Cavanaugh and Neath (1999).

For each $m \in \{1, \cdots, M\}$, the marginal distribution $\nu_m$
of the response vector $\bY$ conditional on model $\mathfrak{M}_m$ has
the density function
\begin{equation} \label{065}
\frac{d \nu_m}{d \mu_0}\left(\bz\right) = \int \exp\left[\ell_n(\bz,
  \bbeta)\right] d\mu_{\mathfrak{M}_m}(\bbeta), \quad \bz \in
\mathbb{R}^n.
\end{equation}
The KL divergence of $\nu_m$ from the true model $G_n$ is given by
\begin{equation} \label{066}
I\left(g_n; \frac{d \nu_m}{d \mu_0}\right) = E \log g_n(\bY) - E \log \int \exp\left[\ell_n(\bY, \bbeta)\right]
d\mu_{\mathfrak{M}_m}(\bbeta),
\end{equation}
which leads to the following proposition.

\begin{proposition} \label{prop1}
\begin{itemize}
\item[a)]
For each $m \in \{1, \cdots, M\}$, we have
\begin{equation} \label{067}
E S(\bY, \mathfrak{M}_m; F_n) = -I\left(g_n; \frac{d \nu_m}{d
  \mu_0}\right) + \log \alpha_{\mathfrak{M}_m} + E \log g_n(\bY),
\end{equation}
where the expectation is taken with respect to the true distribution $G_n$.

\item[b)] Assume $\alpha_{\mathfrak{M}_1} = \cdots = \alpha_{\mathfrak{M}_M}$. Then we have
\[
\arg \max_{m \in \{1, \cdots, M\}} E S(\bY, \mathfrak{M}_m; F_n) =
\arg \min_{m \in \{1, \cdots, M\}} I\left(g_n; \frac{d \nu_m}{d
  \mu_0}\right).
\]
\end{itemize}
\end{proposition}

Proposition \ref{prop1} holds true regardless of whether or not the true model is in
the set of
candidate models. We see that for each $m \in \{1, \cdots, M\}$,
$-S(\by, \mathfrak{M}_m; F_n) + \log \alpha_{\mathfrak{M}_m}$ gives, up to a
common additive constant, an unbiased estimate of $I\left(g_n;
\frac{d \nu_m}{d \mu_0}\right)$. This shows that the Bayesian
principle of model selection can be restated as
choosing the model that minimizes the KL divergence of the marginal
distribution of the response vector $\bY$ from its true distribution, provided
that we assign equal prior probabilities on the $M$ competing models.

Given a set of candidate models $\left\{\mathfrak{M}_m: m = 1,
\cdots, M\right\}$ of form (\ref{003}), we can  also
compare them using the KL
divergence principle of model selection introduced by Akaike (1973,
1974), who studied the case of i.i.d. observations with correctly
specified models.  With the  sequence of QMLE's
$\{\hbbeta_{n, m}: m = 1, \cdots, M\}$ based on these models, Akaike's principle of model selection is to choose
the model $\mathfrak{M}_{m_0}$ that minimizes the KL
divergence $I(g_n; f_n(\cdot, \hbbeta_{n, m}))$ of the fitted
model  $F_n(\cdot, \hbbeta_{n, m})$ from the true model $G_n$, that is,
$
m_0 = \arg \min_{m \in \{1, \cdots, M\}} I(g_n; f_n(\cdot, \hbbeta_{n, m})).
$
This principle is equivalent to choosing the
model  $\mathfrak{M}_{m_0}$ that maximizes the expected
log-likelihood with the expectation taken with respect to an
independent copy of $\bY$. For notational simplicity, we drop the
subscript $m$ hereafter when there is no confusion. Using the
asymptotic theory of MLE, Akaike (1973)
showed for the case of i.i.d. observations that
the expected log-likelihood can be asymptotically
expanded as $\ell_n(\by, \hbbeta_n) -
|\mathfrak{M}|$, which leads to the seminal AIC for comparing
models:
\begin{equation} \label{062}
\mathrm{AIC}(\by, \mathfrak{M}) = -2 \ell_n(\by, \hbbeta_n) + 2 |\mathfrak{M}|.
\end{equation}

\section{Asymptotic expansions of model selection principles in
  misspecified models} \label{sec3}
We here derive a novel asymptotic expansion of the
Bayesian principle in  misspecified GLMs. Our derivation indicates
that, in contrast to traditional AIC and BIC, directly incorporating
the effect of model misspecification in a
model selection criterion can be beneficial. The KL divergence
principle has been considered by many researchers in
the setting of model misspecification (see Section
\ref{sec3.3}); the derivations were, however, mostly for the case
of i.i.d. observations and often by heuristic arguments. We also give a rigorous
asymptotic expansion of this principle in misspecified
GLMs with deterministic design matrices.

\subsection{Generalized BIC in misspecified models} \label{sec3.1}
Let us consider the log-marginal likelihood
$S(\by, \mathfrak{M}; F_n)$, as defined in (\ref{083}), in
misspecified models. In GLM (\ref{003}), the
parameter vector $\btheta \in
\mathbb{R}^n$ for the response vector $\bY$ is given by $\btheta = \bX
\bbeta$. Under this working model and some regularity conditions, as in
Schwartz (1978), we can apply the classical Laplace approximation
(that is, Taylor-expanding the exponent
to the second order and using the Gaussian integration to obtain
necessary terms) to show that with probability tending to one,
\begin{align}\label{bic}
\nonumber
S(\by, \mathfrak{M}; F_n) & = \ell_n(\by, \hbbeta_n)+ \log \int
\exp\left[-(\bbeta-\hbbeta_n)\t \bA_n(\hbbeta_n) (\bbeta-\hbbeta_n)/2\right] d\bbeta + C_n' \\
& = \ell_n(\by, \hbbeta_n) - \frac{\log
  n}{2} |\mathfrak{M}| - \frac{1}{2} \log |n^{-1} \bA_n(\hbbeta_n)| + C_n,
\end{align}
where $\hbbeta_n$ is the maximizer of $\ell_n(\by,\bbeta)$,
$\bA_n(\bbeta) = -\partial^2 \ell_n(\by, \bbeta)/\partial \bbeta^2
= \bX\t \Sig(\bX \bbeta) \bX$, and $C_n'$ and $C_n$ are both bounded quantities. When the model is misspecified, it is, however, no longer clear whether the $C_n$
term is still bounded. Moreover, the second-order term $-\frac{1}{2} \log |n^{-1} \bA_n(\hbbeta_n)|$ in
(\ref{bic}) depends on the scale of the design matrix $\bX$ and thus is not very meaningful. We will derive an asymptotic expansion of the Bayes factor in misspecified models using the asymptotic theory of the QMLE presented in Section \ref{sec4.2}.

Let $\mu_\mathfrak{M}$ be the prior distribution on $\bbeta \in
\mathbb{R}^d$ for model $\mathfrak{M}$ and $\pi(\bbeta) =
\frac{d\mu_\mathfrak{M}}{d\mu_0} (\bbeta)$ its prior density relative
to the Lebesgue measure $\mu_0$ on $\mathbb{R}^d$. In high-dimensional Bayesian modeling, we need to put some basic
constraint on the prior distribution in order to prevent it from
assigning a diminishingly low probability in a vicinity of the
theoretically best model under the KL divergence. Let $\delta_n$ be a
sequence of diverging positive numbers. By the proof of consistency of the QMLE
$\hbbeta_n$ in Theorem \ref{thm5}, we have the concentration probability bound
\[ P\{\hbbeta_n \in N_n(\delta_n)\} \geq 1 - \frac{4 d}{c^2 \delta_n^2}, \]
where $N_n(\delta_n) = \{\bbeta \in \mathbb{R}^d: \|\bB_n^{1/2} (\bbeta - \bbeta_{n, 0})\| \leq \delta_n\}$ is an asymptotically shrinking neighborhood of the parameter vector $\bbeta_{n, 0}$ of the best GLM $F_n(\cdot, \bbeta_{n, 0})$ and $c$ is some positive constant specified in Condition \ref{cond2} in Section \ref{sec4.1}. We denote by
\[ \widetilde{N}_n(\delta_n) = \{\bbeta
\in \mathbb{R}^d: \|\bB_n^{1/2} (\bbeta - \hbbeta_n)\| \leq
\delta_n\} \]
a neighborhood of the QMLE $\hbbeta_n$ and $\widetilde{N}_n^c(\delta_n)$ the
complement of $\widetilde{N}_n(\delta_n)$. Clearly
$\widetilde{N}_n(\delta_n) \subset N_n(2 \delta_n)$ conditional on the above event $\{\hbbeta_n \in N_n(\delta_n)\}$. The key idea of
deriving the asymptotic expansion of the log-marginal likelihood
$S(\by, \mathfrak{M}; F_n)$ is to expand the quasi-log-likelihood
function $\ell_n(\by, \bbeta)$ at the QMLE $\hbbeta_n$ and show that
the integral with the prior distribution on the neighborhood
$\widetilde{N}_n(\delta_n)$ dominates that on
$\widetilde{N}_n^c(\delta_n)$. Note that the neighborhood $N_n(2 \delta_n)$ of $\bbeta_{n, 0}$ shrinks
asymptotically with the sample size $n$, meaning that the constraint
on the prior becomes weaker asymptotically.

For any $\bbeta \in \mathbb{R}^d$, we define a quantity
\begin{equation} \label{084}
\ell^*_n(\by, \bbeta) = \ell_n(\by, \bbeta)
- \ell_n(\by, \hbbeta_n),
\end{equation}
which is the deviation of the quasi-log-likelihood from its
maximum. It follows from (\ref{083}) and (\ref{084}) that
\begin{equation} \label{023}
S(\by, \mathfrak{M}; F_n) =  \ell_n(\by,
\hbbeta_n) + \log E_{\mu_\mathfrak{M}} [U_n(\bbeta)^n] + \log \alpha_\mathfrak{M},
\end{equation}
where $U_n(\bbeta) = \exp\left[n^{-1} \ell^*_n(\by,
  \bbeta)\right]$. We present the asymptotic expansion of the above term $\log
E_{\mu_\mathfrak{M}} [U_n(\bbeta)^n]$ in the following theorem.

\begin{theorem} \label{thm1}
Under Conditions \ref{cond1}, \ref{cond2}, and \ref{cond5} in Section
\ref{sec4.1}, we have with probability tending to one,
\begin{equation} \label{034}
S(\by, \mathfrak{M}; F_n) = \ell_n(\by, \hbbeta_n) -\frac{\log n}{2}
|\mathfrak{M}| + \frac{1}{2} \log |\bH_n| + \log \alpha_\mathfrak{M} + \frac{\log (2 \pi)}{2} |\mathfrak{M}| +\log c_n + o(1),
\end{equation}
where $\bH_n = \bA_n^{-1} \bB_n$ and $c_n \in [c_1, c_2]$.
\end{theorem}

Following the asymptotic expansion in the above theorem, we introduce the generalized BIC (GBIC) as follows.

\begin{definition}
We define GBIC of the competing model $\mathfrak{M}$ by
\begin{equation} \label{035}
\mathrm{GBIC}(\by, \mathfrak{M}; F_n) = -2 \ell_n(\by, \hbbeta_n) + \left(\log
  n\right) |\mathfrak{M}| - \log |\widehat{\bH}_n|,
\end{equation}
where $\widehat{\bH}_n$ is an estimate of the covariance contrast matrix $\bH_n$ given in Section \ref{sec3.4}.
\end{definition}

Compared with BIC, GBIC takes into account model misspecification
explicitly. The second term in GBIC is the same as that in BIC and
penalizes the model complexity. The third term $-\log |\widehat{\bH}_n|$ in GBIC is, however, not always nonnegative and thus is not necessarily a penalty term on model misspecification. When the model is correctly specified, we would
expect $\widehat{\bH}_n \approx I_d$ since $\bH_n = I_d$ due to $\bA_n = \bB_n$
and thus $\log |\widehat{\bH}_n| \approx \log |I_d| = 0$. In this
sense, GBIC introduced above generalizes BIC. We will see in Section \ref{sec5} that the second-order term in GBIC can lead to improved performance in small samples even in correctly specified models.

\subsection{$\mbox{GBIC}_p$ in misspecified models} \label{sec3.2}
Assigning a reasonable prior probability to each candidate model is an
art. One usually chooses to give a higher prior probability to a model
with lower complexity.
In this section, we consider enlightening this prior choice
with the KL divergence principle. For each
candidate model $\mathfrak{M}_m$ with  $1 \leq m \leq M$, denote by
\begin{equation} \label{006}
D_m = E \left[I(g_n; f_n(\cdot, \hbbeta_{n, m})) - I(g_n; f_n(\cdot, \bbeta_{n, m, 0}))\right]
\end{equation}
the expected excess KL divergence of the best fitted model
$F_n(\cdot, \hbbeta_{n, m})$ from the true model $G_n$ relative to
that of the theoretically best model $F_n(\cdot, \bbeta_{n, m, 0})$
from the true model. It follows from the definition of $\bbeta_{n, m,
  0}$ in Theorem \ref{thm4} that
\[ I(g_n; f_n(\cdot, \hbbeta_{n, m})) - I(g_n; f_n(\cdot, \bbeta_{n,
  m, 0})) \geq 0, \]
which entails that $D_m$  is  nonnegative. It is of interest to consider
the following prior assignment: $\alpha_{\mathfrak{M}_m} \propto
e^{-D_m}$, which prefers models with less variability. In the case when
the model family is correctly specified, this prior probability is asymptotically
proportional to $e^{-d/2}$, where $d$ is the number of covariates in the
model.  By the proof of
Theorem \ref{thm3}, we have $\log \alpha_\mathfrak{M} = -\frac{1}{2}
\tr(\bH_n) + o(1)$ up to a common additive constant, which along with Theorem \ref{thm1} leads to the
following asymptotic expansion of the Bayes factor.

\noindent
\begin{theorem} \label{thm2}
Under Conditions \ref{cond1}--\ref{cond6} in Section \ref{sec4.1}, if
the prior probabilities $\alpha_{\mathfrak{M}_m} = C e^{-D_m}$ with $C$ some normalization constant, we
have with probability tending to one,
\begin{equation} \label{034}
S(\by, \mathfrak{M}; F_n) = \ell_n(\by, \hbbeta_n) -\frac{\log n}{2}
|\mathfrak{M}| -\frac{1}{2} \tr(\bH_n) + \frac{1}{2} \log |\bH_n| + \frac{\log (2 \pi)}{2} |\mathfrak{M}| +\log (c_n C) + o(1),
\end{equation}
where $\bH_n = \bA_n^{-1} \bB_n$ and $c_n \in [c_1, c_2]$.
\end{theorem}

Following the asymptotic expansion in the above theorem, we introduce the generalized BIC with prior probability ($\mbox{GBIC}_p$) as follows.

\begin{definition}
We define $\mbox{GBIC}_p$ of the competing model $\mathfrak{M}$ by
\begin{equation} \label{004}
\mbox{GBIC}_p(\by, \mathfrak{M}; F_n) = -2 \ell_n(\by, \hbbeta_n) +\left(\log
  n\right) |\mathfrak{M}| + \tr(\widehat{\bH}_n) - \log |\widehat{\bH}_n|,
\end{equation}
where $\widehat{\bH}_n$ is an estimate of the covariance contrast matrix $\bH_n$ given in Section \ref{sec3.4}.
\end{definition}

Let us consider the case of correctly specified models. As discussed in Section \ref{sec3.4}, we can construct a consistent estimate $\widehat{\bH}_n$ of the covariance contrast matrix $\bH_n = I_d$, which entails that the second-order term $\tr(\widehat{\bH}_n) - \log |\widehat{\bH}_n|$ is asymptotically close to $|\mathfrak{M}|$. This term is asymptotically dominated by the penalty term $(\log n) |\mathfrak{M}|$ on model complexity. In view of (\ref{004}), this observation suggests that under the uniform asymptotic dominance of the second-order term, the model selection criterion $\mbox{GBIC}_p$ is asymptotically equivalent to BIC in correctly specified models, leading to model selection consistency; similar result also holds for GBIC. The $\mbox{GBIC}_p$ also enjoys a representation shown in the following proposition.

\begin{proposition} \label{prop2}
When $\widehat{\bH}_n$ in (\ref{004}) is chosen to be $\widehat{\bA}_n^{-1} \widehat{\bB}_n$ for some $d \times d$ symmetric positive definite matrices $\widehat{\bA}_n$ and $\widehat{\bB}_n$, we have
\begin{equation} \label{071}
\mbox{GBIC}_p(\by, \mathfrak{M}; F_n) = - 2 \ell_n(\by, \hbbeta_n) + (1 +
  \log n) |\mathfrak{M}| + 2 I \left[N(\bzero, \widehat{\bB}_n);
  N(\bzero, \widehat{\bA}_n)\right],
\end{equation}
where $I [N(\bzero, \widehat{\bB}_n); N(\bzero, \widehat{\bA}_n)] =
\frac{1}{2} [\tr(\widehat{\bH}_n) - \log |\widehat{\bH}_n| - d]$ is
the KL divergence of the $d$-variate Gaussian distribution $N(\bzero,
\widehat{\bA}_n)$ from the $d$-variate Gaussian distribution
$N(\bzero, \widehat{\bB}_n)$.
\end{proposition}

The above representation shows that $\mbox{GBIC}_p$ in misspecified models admits a natural decomposition of the
form
\begin{equation} \label{098}
\text{goodness of fit + model complexity + model misspecification},
\end{equation}
where the first term is the negative maximum quasi-log-likelihood and the two penalty terms on model
complexity and model misspecification are both nonnegative. However,
neither the plain GBIC nor GAIC, which will be introduced in Section
\ref{sec3.3}, has such a decomposition. The model misspecification
term in $\mbox{GBIC}_p$ has a natural interpretation. Under the true
model $G_n$, we have $E (\bX\t \bY) = \bX\t E \bY$ and $\cov(\bX\t \bY) = \bX\t
\cov(\bY) \bX = \bB_n$. Under some regularity conditions, an application of the central limit theorem shows that the asymptotic distribution of
$\bX\t (\bY-\bmu)$ with $\bmu = E_{G_n} \bY$ is close to $N(\bzero, \bB_n)$. On the other hand, if we misspecify the true model as the GLM
$F_n(\cdot, \bbeta)$, the QMLE $\hbbeta_n$ asymptotically converges to $\bbeta_{n,0}$
such that $F_n = F_n(\cdot,\bbeta_{n,0})$ is the best model in the posited
GLM family to approximate the true model $G_n$ under the
KL divergence.
It follows from (\ref{011}) and (\ref{073}) that
$E_{F_n} (\bX\t \bY) = \bX\t E_{F_n} \bY = \bX\t \bmu$ and $\cov_{F_n}(\bX\t \bY) = \bX\t
\Sig(\bX \bbeta_{n, 0}) \bX = \bA_n$. Under some regularity conditions, the
asymptotic distribution of $\bX\t (\bY-\bmu)$ is now close to
$N(\bzero, \bA_n)$. Thus the KL divergence
of $N(\bzero, \bA_n)$ from $N(\bzero,
\bB_n)$ provides a natural measure of model misspecification. When the
model is correctly specified, we would expect $\widehat{\bA}_n \approx
\widehat{\bB}_n$ since $\bA_n = \bB_n$, and thus $I [N(\bzero, \widehat{\bB}_n); N(\bzero,   \widehat{\bA}_n)] \approx 0$. We will show in Section \ref{sec5} that the second-order term in $\mbox{GBIC}_p$ can also lead to improved performance in small samples even in correctly specified models.

\subsection{Generalized AIC in misspecified models} \label{sec3.3}
We now consider the KL divergence principle in misspecified models. AIC has been studied by many researchers, mainly for the case
of i.i.d observations. Takeuchi (1976) generalized the derivation of
AIC and obtained a term $\tr(\bA^{-1} \bB)$ in place of the dimension
of the model, where $\bB$ and $\bA$ are two  matrices that involve the
first and second derivatives of the log-likelihood function,
respectively. This term was also obtained by Stone (1977) in deriving
the asymptotic expansion for the cross-validation. In fact,
$\tr(\bA^{-1} \bB)$ is the well-known Lagrange-multiplier test
statistic. See also Hosking (1980), Shibata (1989), Konishi and
Kitagawa (1996), Burnham and Anderson (1998), and Bozdogan (2000).

In view of (\ref{008}), the KL divergence of the fitted model $F_n(\cdot, \hbbeta_{n, m})$ from the true model $G_n$ can be written as
\begin{equation} \label{012}
I(g_n; f_n(\cdot, \hbbeta_{n, m})) = E \log g_n(\widetilde{\bY}) - \eta_n(\hbbeta_{n, m}),
\end{equation}
where $\eta_n(\bbeta) = E \ell_n(\widetilde{\bY}, \bbeta)$ and
$\widetilde{\bY}$ is an independent
copy of $\bY$. Thus the best model chosen by Akaike's principle of model
selection is
$
m_0 = \arg \max_{m \in \{1, \cdots, M\}} \eta_n(\hbbeta_{n, m}).
$
For notational simplicity, we drop the subscript $m$. We also drop the last constant term in the
quasi-log-likelihood $\ell_n(\by, \bbeta)$ defined in (\ref{002}), which
does not depend on $\bbeta$, and redefine it as
$
\ell_n(\by, \bbeta) = \by\t \bX \bbeta - \bone\t \bb(\bX \bbeta).
$
In view of (\ref{012}), for the purpose of model comparison, we need to estimate the quantity of the expected log-likelihood $\eta_n(\hbbeta_n)$. In practice, we have only
a single sample $\by$ available. Clearly the maximum quasi-log-likelihood $\ell_n(\by, \hbbeta_n)$ tends
to overestimate $\eta_n(\hbbeta_n)$ since we would use the same data
twice, that is, in both estimation and model evaluation. We present the asymptotic expansion of $E \eta_n(\hbbeta_n)$ in the following theorem.

\begin{theorem} \label{thm3}
Under Conditions \ref{cond1}--\ref{cond4} and \ref{cond6} in Section \ref{sec4.1}, we have
\begin{equation} \label{001}
E \eta_n(\hbbeta_n) = E \ell_n(\by, \hbbeta_n) -
\tr\left(\bH_n\right) + o(1),
\end{equation}
where $\bH_n = \bA_n^{-1} \bB_n$ and both expectations are taken with respect to the true distribution $G_n$.
\end{theorem}

As mentioned before, a correction term of form $\tr\left(\bA_n^{-1}
\bB_n\right)$ is well known in the literature, but has been derived usually for the case of i.i.d. observations and often by heuristic
arguments. Following the asymptotic expansion in the above theorem, we
have the generalized AIC (GAIC) as follows.

\begin{definition} The GAIC of the competing model $\mathfrak{M}$ is
\begin{equation} \label{017}
\mathrm{GAIC}(\by, \mathfrak{M}; F_n) = -2 \ell_n(\by, \hbbeta_n) +
2 \tr\left(\widehat{\bH}_n\right),
\end{equation}
where $\widehat{\bH}_n$ is an estimate of the covariance contrast matrix $\bH_n$ given in Section \ref{sec3.4}.
\end{definition}

The GAIC shares some common features with the GBIC. It incorporates the
effects of model complexity and
model misspecification in a single term. When $\widehat{\bH}_n$ in (\ref{017}) is chosen to be $\widehat{\bA}_n^{-1} \widehat{\bB}_n$ for some $d \times d$ symmetric positive definite matrices $\widehat{\bA}_n$ and $\widehat{\bB}_n$, we have $\tr(\widehat{\bH}_n) = \tr(\widehat{\bB}_n^{1/2}
\widehat{\bA}_n^{-1} \widehat{\bB}_n^{1/2}) \geq 0$, showing that the term added to the negative maximum
quasi-log-likelihood in GAIC is indeed
a penalty term. When the model is correctly specified, it holds that $\tr(\bH_n) = \tr(I_d) = d = |\mathfrak{M}|$ since $\bH_n = I_d$ due to $\bA_n = \bB_n$, which is the number
of predictors in the GLM (\ref{003}). Thus we would expect $\tr(\widehat{\bH}_n) \approx \tr\left(I_d\right) =
|\mathfrak{M}|$, which reflects the penalty on model complexity as in AIC. As mentioned before, many extensions of AIC under different settings share the same form as in (\ref{017}). We will also see in Section \ref{sec5} that the second-order term in GAIC can lead to improved performance in small samples even in correctly specified models. We point out that equation (4.5) in Stone's derivation (Stone, 1977) suggests that under the same working models, the GAIC and cross-validation are asymptotically equivalent under some weak conditions, where the cross-validation uses the log-density assessment.

\subsection{Estimation of covariance contrast matrix $\bH_n$} \label{sec3.4}
The model selection criteria GAIC, GBIC, and $\mbox{GBIC}_p$ require an estimate $\widehat{\bH}_n$ of the covariance contrast matrix $\bH_n$, a contrast between the covariance structures in the misspecified model and in the true model. Specifically, we need to estimate the trace $\tr\left(\bH_n\right)$ for GAIC, the log-determinant $\log |\bH_n|$ for GBIC, and both quantities for $\mbox{GBIC}_p$. Estimating these quantities in misspecified models is nontrivial. In this paper, we consider two specific estimators for them and illustrate the performance of the resulting model selection criteria with extensive numerical studies. More accurate estimation of $\tr\left(\bH_n\right)$ and $\log |\bH_n|$, especially in high dimensions, demands rigorous theoretical studies, which we will not address in the paper. Through our numerical studies, we found that estimating the trace of the covariance contrast matrix $\bH_n$ and estimating its log-determinant exhibit different behavior. To show the utility of the new model selection criteria, we estimate the two quantities separately with two different estimators $\widehat{\bH}_n$: a simple estimator for $\log |\bH_n|$ and a bootstrap estimator for $\tr\left(\bH_n\right)$.

To estimate $\log |\bH_n|$, we consider the estimate $\widehat{\bH}_n$ of the plug-in form $\widehat{\bA}_n^{-1} \widehat{\bB}_n$ with $\widehat{\bA}_n$ and $\widehat{\bB}_n$ some $d \times d$ symmetric positive definite matrices. Since the QMLE $\hbbeta_n$ provides a consistent estimator of $\bbeta_{n, 0}$ in the best misspecified GLM $F_n(\cdot, \bbeta_{n, 0})$, a natural estimate of matrix $\bA_n$ is given by
\begin{equation} \label{087}
\widehat{\bA}_n = \bA_n(\hbbeta_n) = \bX\t \Sig(\bX \hbbeta_n) \bX.
\end{equation}
When the model is correctly specified or more generally $E\bY = \bmu(\bX \bbeta_{n,0})$, the following simple estimator
\begin{equation} \label{088}
\widehat{\bB}_n = \bX\t \diag\left\{\left[\by - \bmu(\bX
  \hbbeta_n)\right] \circ \left[\by - \bmu(\bX
  \hbbeta_n)\right]\right\} \bX,
\end{equation}
with ``$\circ$'' representing the Hadamard (componentwise) product, provides an asymptotically unbiased estimator of $\bB_n$. In view of (\ref{089}), we adopt the form of the diagonal matrix in (\ref{088}) for estimating $\log |\bH_n|$. This form seems effective in ensuring the positive definiteness and well-conditionedness of the matrix, which are crucial for effectively estimating the log-determinant. We refer to the resulting estimator $\widehat{\bH}_n$ as the simple estimator. The formulas for matrices $\widehat{\bA}_n$ and $\widehat{\bB}_n$ in (\ref{087})--(\ref{088}) are provided in Appendix \ref{secB}, when the working model $F_n$ is chosen to be one of the three commonly used GLMs: the linear regression model, logistic regression model, and Poisson regression model. The following theorem justifies the consistency of the above simple estimator $\widehat{\bH}_n = \widehat{\bA}_n^{-1} \widehat{\bB}_n$ when the model is correctly specified.

\begin{theorem} \label{thm7}
Assume that the model is correctly specified, Conditions \ref{cond1}, \ref{cond2}, and \ref{cond7} in Section
\ref{sec4.1} hold, and the smallest and largest eigenvalues of $n^{-1} \bB_n$ are bounded away from 0 and $\infty$. Then $\widehat{\bH}_n = \widehat{\bA}_n^{-1} \widehat{\bB}_n = I_d + o_P(1)$, where $\widehat{\bA}_n$ and $\widehat{\bB}_n$ are given in (\ref{087}) and (\ref{088}), respectively.
\end{theorem}

In Theorem \ref{thm7}, the assumption that the smallest and largest eigenvalues of $n^{-1} \bB_n$ are bounded away from 0 and $\infty$ can be relaxed to converge to 0 or diverge to $\infty$, which will affect the consistency rate. Although it is positive definite, the simple estimator $\widehat{\bB}_n$ in (\ref{088}) may not be most suitable for the trace estimation. The general bias in such an estimator can affect more the performance of GAIC than that of GBIC and $\mbox{GBIC}_p$, due to the presence of the penalty term $\left(\log n\right) d$ on model complexity in the latter two information criteria. For example, $\tr\left\{(1 + c) \bH_n\right\} - \tr\left(\bH_n\right) = c \tr\left(\bH_n\right)$ is of the same order as $\tr\left(\bH_n\right)$ for some nonzero constant $c$ close to zero, while $\log |(1 + c) \bH_n| - \log |\bH_n| = [\log(1 + c)] d$ can be asymptotically dominated by the term $\left(\log n\right) d$. To mitigate the bias in estimating $\bB_n$, it is sensible to construct the residual vector with a large model. This approach generally gives rise to a less biased estimator for $\bB_n$. Motivated by these considerations, to estimate $\tr\left(\bH_n\right)$ we introduce the estimate $\widehat{\bH}_n = \widehat{\bA}_n^{-1} \widehat{\bB}_n$ with $\widehat{\bA}_n$ given above and $\widehat{\bB}_n$ constructed using the following bootstrap procedure (Efron, 1979). To reduce the model bias, we first use a large model, with regularization when needed, and the idea of cross-validation to construct a proxy $\bd$ of the deviation vector $\bY - E \bY$. Then, we draw some bootstrap samples, each of size $n$, from the $n$ observed vectors of the covariates and constructed deviations. For each bootstrap sample $(\bX^*, \bd^*)$, we calculate the covariance vector $(\bX^*)\t \bd^*$. Finally, we obtain the bootstrap estimate $\widehat{\bB}_n$ as the covariance matrix of the covariance vectors constructed from the bootstrap samples. We refer to the resulting estimator $\widehat{\bH}_n$ as the bootstrap estimator.

The estimates for $\tr\left(\bH_n\right)$ and $\log |\bH_n|$ introduced above may not be the most effective ones. In practice, one can construct other estimates of both quantities by, for example, treating the squared residual $[\by - \bmu(\bX \hbbeta_n)] \circ [\by - \bmu(\bX \hbbeta_n)]$ as a function of the corresponding fitted value  $\bmu(\bX \hbbeta_n)$ and some covariates, and doing some local smoothing. Another possible way of estimating $\bH_n$ is using the bootstrap as suggested by Theorem \ref{thm6}. We show in Theorem \ref{thm6} that under certain regularity conditions, the QMLE $\hbbeta_n$ is asymptotically normal with mean $\bbeta_{n, 0}$ and covariance matrix $\bC_n^{-1} (\bC_n^{-1})\t = \bA_n^{-1} \bB_n \bA_n^{-1}$, where $\bC_n = \bB_n^{-1/2} \bA_n$. In view of $\bH_n = (\bA_n^{-1} \bB_n \bA_n^{-1}) \bA_n$, this asymptotic normality motivates us to construct the estimate $\widehat{\bH}_n$ as $\widehat{\text{Cov}}(\hbbeta_n) \widehat{\bA}_n$ using the bootstrap estimate $\widehat{\text{Cov}}(\hbbeta_n)$ of the covariance matrix of the QMLE $\hbbeta_n$. Regularization methods for large covariance matrix estimation can also be exploited to improve the estimation of both $\tr\left(\bH_n\right)$ and $\log |\bH_n|$.

\section{Technical conditions and results} \label{sec4}

\subsection{Technical conditions} \label{sec4.1}

Denote by $\lambda_{\min}(\cdot)$ and  $\lambda_{\max}(\cdot)$ the smallest and largest eigenvalues of a given matrix, respectively, and $\|\cdot\|$ the Euclidean norm of a vector as well as the matrix operator norm. We need the following mild regularity conditions.

\begin{assumption} \label{cond1}
$b(\theta)$ is twice differentiable with
$b''(\theta)$ always positive and $\bX$ has full
column rank $d$.
\end{assumption}

\begin{assumption} \label{cond2}
$\lambda_{\min}(\bB_n) \rightarrow \infty$ as $n \rightarrow \infty$ and there exists some $c > 0$ such that for any $\delta > 0$, $\min_{\bbeta \in
N_n(\delta)} \lambda_{\min}[\bV_n(\bbeta)] \geq c$ for all sufficiently large $n$, where $\bV_n(\bbeta) = \bB_n^{-1/2} \bA_n(\bbeta)
\bB_n^{-1/2}$ and $N_n(\delta) = \{\bbeta \in \mathbb{R}^d: \|\bB_n^{1/2} (\bbeta - \bbeta_{n, 0})\| \leq \delta\}$.
\end{assumption}

\begin{assumption} \label{cond3}
For any $\delta > 0$, $\max_{\bbeta_1, \cdots, \bbeta_d \in
  N_n(\delta)} \|\widetilde{\bV}_n(\bbeta_1, \cdots, \bbeta_d) -
\bV_n\| = o(1)$, where $\bV_n =
\bV_n(\bbeta_{n, 0}) = \bB_n^{-1/2} \bA_n \bB_n^{-1/2}$ and $\widetilde{\bV}_n(\bbeta_1, \cdots, \bbeta_d) =
\bB_n^{-1/2} \widetilde{\bA}_n(\bbeta_1, \cdots, \bbeta_d)
\bB_n^{-1/2}$ with the $j$-th row of $\widetilde{\bA}_n(\bbeta_1, \cdots,
\bbeta_d)$ the
corresponding row of $\bA_n(\bbeta_j)$ for each $j = 1, \cdots, d$.
\end{assumption}

\begin{assumption} \label{cond4}
$\max_{i = 1}^n E |Y_i - E Y_i|^3 = O(1)$ and $
  \sum_{i = 1}^n (\bx_i\t \bB_n^{-1} \bx_i)^{3/2} = o(1)$.
\end{assumption}

\begin{assumption} \label{cond5}
The prior density $\pi(\bbeta)$ satisfies
\[ \inf_{\bbeta \in N_n(2 \delta_n)} \pi\left(\bbeta\right) \geq c_1
\ \text{ and } \ \sup_{\bbeta
  \in \mathbb{R}^d} \pi\left(\bbeta\right) \leq c_2 \]
for some positive constants $c_1, c_2$ and a sequence of positive
numbers $\delta_n$ satisfying $\delta_n^{-1} = o(\log^{-1/2} n)$ and
$\max_{\bbeta \in N_n(2 \delta_n)} \max\{|\lambda_{\min}(\bV_n(\bbeta)
- \bV_n)|, |\lambda_{\max}(\bV_n(\bbeta) - \bV_n)|\} = o(1)$, and $\lambda_{\max}(\bV_n)$ is of a polynomial order of $n$.
\end{assumption}

\begin{assumption} \label{cond6}
$\sup_{\bdelta} \max_{j, k, l} |\partial^3 \eta_n(\bbeta_{n, 0} + n^{1/2} \bC_n^{-1} \bdelta)/\partial \delta_j \partial \delta_k \partial \delta_l| = o(n^{3/2})$ and $E \|\bC_n (\hbbeta_n - \bbeta_{n, 0})\|^{3 + \alpha} = O(1)$ for some $\alpha > 0$, where $\bC_n = \bB_n^{-1/2} \bA_n$ and $\bdelta = (\delta_1, \cdots, \delta_d)\t$.
\end{assumption}

\begin{assumption} \label{cond7}
There exists some $L > 0$ such that for any $\delta > 0$, the functions $n^{-1} \bA_n(\bbeta)$, $n^{-1} \bX\t \diag\{|\bmu(\bX \bbeta) - \bmu(\bX \bbeta_{n, 0})|\} \bX$, and $n^{-1} \bX\t \diag\{[\bmu(\bX \bbeta) - \bmu(\bX \bbeta_{n, 0})] \circ [\bmu(\bX \bbeta) - \bmu(\bX \bbeta_{n, 0})]\} \bX$ are Lipschitz with a Lipschitz constant $L$ in the neighborhood $N_n(\delta)$. In addition, $n^{-1/2} \sum_{i = 1}^n \{\bx_i \bx_i\t (Y_i - E Y_i)^2\}$ has bounded second central moment and $P(\max_{i = 1}^n |Y_i - E Y_i| > K_n) \rightarrow 0$ as $n \rightarrow \infty$ for some $K_n = o(n^{1/2})$.
\end{assumption}

Condition \ref{cond1} contains some basic smoothness and regularity conditions to ensure the uniqueness of the QMLE $\hbbeta_n$ by a standard strict concavity argument and the uniqueness of the best misspecified GLM $F_n(\cdot, \bbeta_{n, 0})$ under the KL divergence in Theorem \ref{thm4}.

Condition \ref{cond2} is for establishing the consistency of the QMLE $\hbbeta_n$ in Theorem \ref{thm5}. The first part of Condition \ref{cond2} requires that the smallest eigenvalue of $\bX\t \cov(\bY) \bX$ diverges as $n$ increases. When $\var(Y_i)$'s are bounded away from 0 and $\infty$, this assumption amounts to the usual assumption on the design matrix $\bX$ in linear regression for consistency, that is, $\lambda_{\min}(\bX\t\bX) \rightarrow \infty$ as $n \rightarrow \infty$. Intuitively, the second part of Condition \ref{cond2} means that the covariance structures given by the misspecified GLMs in a shrinking
neighborhood of the best working model $F_n(\cdot, \bbeta_{n, 0})$ cannot be too far away from the covariance
structure under the true model $G_n$. This requirement is analogous to that in importance sampling, that is, the support of the sampling distribution has to cover that of the target distribution. When the model is correctly specified, we have $\bV_n = \bV_n(\bbeta_{n, 0}) = I_d$ since $\bA_n = \bB_n$ by the equivalence of the Hessian and outer product forms for the Fisher information matrix.

Conditions \ref{cond3} and \ref{cond4} facilitate the derivation of the asymptotic normality of the QMLE $\hbbeta_n$ in Theorem \ref{thm6}, where the former is on the continuity of the matrix-valued function $\widetilde{\bV}_n(\bbeta_1, \cdots, \bbeta_d)$ in a shrinking neighborhood of $\bbeta_{n, 0}$ and the latter is a typical moment condition for proving asymptotic normality using Lyapunov's theorem.

The first part of Condition \ref{cond5} is a mild condition on the prior distribution on $\bbeta$ requiring that its prior density is bounded and locally bounded away from zero in a shrinking neighborhood of $\bbeta_{n, 0}$, and facilitates the derivation of the asymptotic expansion of the Bayes factor in misspecified models in Theorems \ref{thm1} and \ref{thm2}. See Section \ref{sec3.1} for more discussions on this condition. The second part of Condition \ref{cond5} ensures that the integral $E_{\mu_\mathfrak{M}} [U_n(\bbeta)^n 1_{\widetilde{N}_n^c(\delta_n)}]$ with the prior distribution outside the neighborhood
$\widetilde{N}_n(\delta_n)$ is negligible in approximating the log-marginal likelihood
$S(\by, \mathfrak{M}; F_n)$.

Condition \ref{cond6} contains some regularity conditions for rigorously deriving the asymptotic expansion of the KL divergence principle in misspecified models in Theorem \ref{thm3}. It is also needed for deriving the asymptotic expansion of the log-prior probability in Theorem \ref{thm2}.

Condition \ref{cond7} collects additional regularity conditions to ensure that the simple estimator $\widehat{\bH}_n = \widehat{\bA}_n^{-1} \widehat{\bB}_n$, with $\widehat{\bA}_n$ and $\widehat{\bB}_n$ given in (\ref{087}) and (\ref{088}), respectively, is a consistent estimator of $\bH_n = I_d$ when the model is correctly specified. The Lipschitz conditions on the three matrix-valued functions in the asymptotically shrinking neighborhood $N_n(\delta)$ are mild and reasonable. The moment and tail conditions on the response variable $Y_i$ are standard.

\subsection{Asymptotic properties of QMLE} \label{sec4.2}
For completeness, we briefly present here some asymptotic properties of the QMLE in misspecified GLMs with deterministic design matrices, which serve as the foundation for
deriving asymptotic expansions of model selection principles. The consistency and asymptotic normality of MLE in correctly specified GLMs were studied in Fahrmeir and Kaufmann (1985). A general theory of maximum likelihood estimation of misspecified models was presented in White (1982), who studied the case of i.i.d. observations from a general distribution and used the KL divergence as a measure of model misspecification. We generalize those results to misspecified GLMs with deterministic design matrices, where the observations may no longer be i.i.d.

\begin{theorem} \label{thm4}
\emph{(Parameter identifiability)}. Under Condition \ref{cond1}, the KL
divergence $I(g_n; f_n(\cdot, \bbeta))$ has a unique global minimum at
$\bbeta_{n, 0} \in \mathbb{R}^d$, which solves the equation
\begin{equation} \label{011}
\bX\t \left[E\bY -
  \bmu(\bX \bbeta)\right] = \bzero.
\end{equation}
\end{theorem}

\begin{theorem} \label{thm5}
\emph{(Consistency)}. Under Conditions \ref{cond1} and \ref{cond2}, the QMLE $\hbbeta_n$
satisfies $\hbbeta_n - \bbeta_{n, 0} = o_P(1)$.
\end{theorem}

\begin{theorem} \label{thm6}
\emph{(Asymptotic normality)}. Under Conditions \ref{cond1}--\ref{cond4}, the QMLE
$\hbbeta_n$ satisfies
$\bC_n (\hbbeta_n - \bbeta_{n, 0}) \toD N(\bzero, I_d)$,
where $\bC_n = \bB_n^{-1/2} \bA_n$.
\end{theorem}

When the model is correctly specified, we have $\bC_n = \bB_n^{-1/2} \bA_n = \bA_n^{1/2}$ since $\bA_n = \bB_n$. Thus in this case, the consistency and asymptotic normality of the QMLE $\hbbeta_n$ in Theorems \ref{thm5} and \ref{thm6} become the conventional asymptotic theory of the MLE. To simplify the technical presentation, the above asymptotic normality is for fixed dimensionality $d$. With more delicate analysis, one can show the asymptotic normality for diverging dimensionality $d$, which is not the focus of the current paper. See, for example, the technical analysis in Fan and Lv (2011) for penalized-maximum likelihood estimation, where the dimensionality is allowed to grow nonpolynomially with the sample size.

\section{Numerical examples} \label{sec5}
In this section, we illustrate the performance of model selection criteria GAIC, GBIC, and $\mbox{GBIC}_p$ in relatively small samples, with the bootstrap estimator for $\tr\left(\bH_n\right)$ and the simple estimator for $\log |\bH_n|$ introduced in Section \ref{sec3.4}. These studies represent some first attempts for obtaining insights into the effect of incorporating model misspecification for model selection; the estimation of $\tr\left(\bH_n\right)$ and $\log |\bH_n|$ still requires rigorous theoretical investigations. Model selection with limited sample size is frequently encountered in practice and the dimension of candidate models can be large due to the large number of covariates. One example in Section \ref{sec5.1.1} is for the case of correctly specified models, and five other examples are for the case of misspecified models. In particular, in Sections \ref{sec5.1.2}-\ref{sec5.1.3} and \ref{sec5.2.3}, we consider high-dimensional linear and logistic regression with interaction, respectively, with Section \ref{sec5.1.3} focusing on the case of moderate effects. Measures of prediction and variable selection are used to evaluate the selected model. We include the classical AIC and BIC as benchmark model selection criteria. With reduced dimensionality given by the set of covariates in the selected model, one can perform more delicate analysis in practice.

\subsection{Linear models} \label{sec5.1}
We consider here three linear models, starting with a correctly specified
model of moderate size and proceeding to high-dimensional sparse linear models with  interaction terms. We observed that GBIC and GBIC$_p$ outperformed other information criteria in all settings, often by a large margin when the sample size is small.

\subsubsection{Best subset linear regression} \label{sec5.1.1}
We simulated 100 data sets from the linear regression model
\begin{equation} \label{053}
\by = \bX \bbeta + \bveps,
\end{equation}
where $\by$ is an $n$-dimensional response vector, $\bX$ is an $n \times p$ design matrix, $\bbeta$ is a $p$-dimensional regression coefficient vector, and the error vector $\bveps \sim N(\bzero,
\sigma^2 I_n)$ is independent of $\bX$. For each simulated data set, the rows of $\bX$ were sampled as i.i.d. copies from $N(\bzero, \Sigma_0)$ with $\Sigma_0 = (0.5^{|i-j|})_{i, j=1,\cdots,p}$. To apply the best subset regression, we chose $p = 6$ and set $\bbeta$ as $\bbeta_0 = (1, -1.25, 0.75, 0, 0, 0)\t$. We considered $n = 20$, $40$, $80$, and $\sigma = 0.25$ and $0.5$, respectively. For each data set, we first applied the best subset regression to build all possible submodels of $\{1, \cdots, 6\}$, each of which numbers represents the corresponding covariate, and then selected the final model using the AIC, BIC, GAIC, GBIC, and $\mbox{GBIC}_p$.

\begin{table}
\caption{Frequency of estimated model size for model in Section \ref{sec5.1.1} over $100$ simulations for different $n$ and $\sigma$; $3^\ast$ corresponds to the true model $\{1, 2, 3\}$}
\centering
\begin{tabular}{clccccccccccccccccc}
\hline
$\sigma$ & Criterion & \multicolumn{5}{c}{Model size ($n = 20$)} & & \multicolumn{5}{c}{Model size ($n = 40$)} & & \multicolumn{5}{c}{Model size ($n = 80$)} \\
\cline{3-7} \cline{9-13} \cline{15-19}
 & & 2 & $3^\ast$ & 4 & 5 & 6 & & 2 & $3^\ast$ & 4 & 5 & 6 & & 2 & $3^\ast$ & 4 & 5 & 6 \\
\cline{3-7} \cline{9-13} \cline{15-19}
0.25 & AIC  & 0  &  55  &  32  &  12   & 1 & & 0  &  55  &  36 &    7  &   2 & & 0  &  71   & 19  &  9   &  1 \\
    & BIC  & 0   & 68  &  24   &  7   &  1 & & 0 &   81  &  15 &    4   &  0 & & 0  &  94   &  5  &   1   &  0 \\
    & GAIC & 0 &   79  &  19   & 1  &   1 & & 0  &  65  &  27  &  7   &  1 & & 0  &  74  &  19  & 6  &   1 \\
    & GBIC & 0  &  78  &  18  &   4 &   0 & & 0   & 83   & 15  &   2 &   0 & & 0  &  95  &   5  &   0  &   0 \\
    & $\mbox{GBIC}_p$  & 0  &  94 &   5  &   1   &  0 & & 0  &  90  &  9  &  1  &   0 & & 0  &  97  &   3  &  0  &   0 \\
\cline{3-7} \cline{9-13} \cline{15-19}
0.5   & AIC  & 0  &  51 &   34  &  10 &    5 & & 0   & 54  &  33  &  11    & 2 & & 0  &  59  &  36   &  3   &  2 \\
    & BIC  & 0   & 65  &  28  &   6 &  1 & & 0 &   80 &  16  &   4   &  0 & & 0 &  85   & 15  &   0  &   0 \\
    & GAIC & 2  &  72   & 21  &  5  &   0 & & 0   & 69  &  21   & 10  &   0 & & 0   & 63  &  34 &    1    & 2 \\
    & GBIC & 0  &  73  &  22   &  5   &  0 & & 0  &  87   & 11 &    2   &  0 & & 0   & 86  &  14  &   0  &   0 \\
    & $\mbox{GBIC}_p$  & 1  &  94  &  5  &  0  &   0 & & 0   & 93   &  7  &  0 &    0 & & 0  &  98  &   2 &    0 &    0 \\
\hline
\end{tabular}
\label{tab5.1.1}
\end{table}

Table 1 summarizes the comparison results of the frequency of estimated model size with $3^\ast$ representing the true underlying sparse model $\supp(\bbeta_0) = \{1, 2, 3\}$ for compactness. The performance of all information criteria generally improves as the sample size increases. AIC tended to select a larger model than the true one even for large sample sizes, while BIC performed well when the sample size is large but also tended to select a larger model than the true one when the sample size is small. GAIC and GBIC improved over AIC and BIC, respectively, across all  sample sizes tested and especially when the sample size is moderate, indicating the effectiveness of the second-order terms even in correctly specified models. When the sample size became larger, GBIC tended to perform better than GAIC. $\mbox{GBIC}_p$ outperformed all other model selection criteria.

\subsubsection{High-dimensional sparse linear regression with interaction} \label{sec5.1.2}
To evaluate how  each information criterion fares in high-dimensional incorrectly specified settings, we simulated  100 data sets from the following linear model:
\begin{equation} \label{092}
\by = \bX \bbeta + \bx_{p + 1} + \bveps,
\end{equation}
where $\bX = (\bx_1, \cdots, \bx_p)$ is an $n \times p$ design matrix, $\bx_{p+1} = \bx_1 \circ \bx_2$ is an interaction term which is the product of the first two covariates, and the rest is the same as in (\ref{053}). To make variable selection easily interpretable in misspecified models, we consider the case that all covariates are independent of each other and thus set $\Sigma_0 = I_p$. We chose $\bbeta$ as $\bbeta_0 = (1, -1.25, 0.75, -0.95, 1.5, 0, \cdots, 0)\t$, $n = 80$, and $\sigma = 0.25$, and considered $p = 50$, $100$, and $200$, respectively. Although the data were generated from model (\ref{092}), we fit the linear regression model (\ref{053}) without interaction. This is an example of misspecified models. In view of (\ref{092}), the true model involves only the first five covariates in a nonlinear form. Since the other covariates are independent of those five covariates, the oracle working model is $\supp(\bbeta_0) = \{1, \cdots, 5\}$, which indicates the linear regression model (\ref{053}) with the first five covariates. Clearly it is unrealistic to implement the best subset regression due to its  computational complexity. Therefore, we used the SICA  regularization method in Lv and Fan (2009) and implemented it with the ICA algorithm in Fan and Lv (2011). For each data set, we first applied the SICA to build a sequence of sparse models and then selected the final model using the same five model selection criteria.

\begin{table}
\caption{Percentiles of prediction error, false positives (FP), and false negatives (FN) as well as selection probability (SP) and inclusion probability (IP) of the SICA estimate for model in Section \ref{sec5.1.2} over $100$ simulations for different $p$}
\centering
\begin{tabular}{clccccccccccc}
\hline
$p$ & Criterion & \multicolumn{3}{c}{Prediction error} & & \multicolumn{3}{c}{FP (FN)} & & SP & & IP \\
\cline{3-5} \cline{7-9}
 & & 25th & 50th & 75th & & 25th & 50th & 75th & &   & &   \\
\cline{3-5} \cline{7-9} \cline{11-11} \cline{13-13}
50 & AIC  & 1.3095  &  1.4433  &  1.6337  & & 4 (0)  &  6 (0)   &  8 (0)   & & 0  & &  1  \\
    & BIC  & 1.1822  &  1.2661  &  1.3924   & & 0 (0)   &  1 (0)  &   4 (0) & & 0.35  & &  0.99   \\
    & GAIC & 1.1843  & 1.2697 &   1.3820   & & 0 (0) &   2 (0) &  4 (0)  & & 0.27  & &  0.99 \\
    & GBIC & 1.1598 &   1.2357   & 1.3401   & & 0 (0) &    1 (0)  &   2 (0)  & & 0.48  & &  0.99 \\
    & $\mbox{GBIC}_p$ & 1.1350  & 1.1873   & 1.2749   & & 0 (0)  &   0 (0)    & 0 (0)  & & 0.79  & &  0.99 \\
    & Oracle & 1.1253  &  1.1600  &  1.2201   & &   0 (0)   &  0 (0)   &  0 (0)  & & 1  & &  1  \\
\cline{3-5} \cline{7-9} \cline{11-11} \cline{13-13}
100 & AIC  & 1.4579  &  1.6807   & 1.9342  & & 7 (0) & 10 (0) & 13 (0)   & & 0  & &  1 \\
    & BIC  & 1.3123 &  1.5229 &   1.7922   & & 2.5 (0) & 7 (0) & 10.5 (0) & & 0.11  & &  1\\
    & GAIC & 1.2038  &  1.2733  &  1.4182  & & 0 (0) &  1 (0) &  2.5 (0) & & 0.37  & &  1 \\
    & GBIC & 1.2193  &  1.3731   & 1.6458   & & 1 (0) & 3 (0) &  7 (0) & & 0.23  & &  1 \\
    & $\mbox{GBIC}_p$ & 1.1253 &   1.2076  &  1.2745   & & 0 (0)  &  0 (0) &  1 (0) & & 0.69  & &  1 \\
    & Oracle & 1.1105  &  1.1639  & 1.2318   & &   0 (0)   &  0 (0)   &  0 (0)  & & 1  & &  1 \\
\cline{3-5} \cline{7-9} \cline{11-11} \cline{13-13}
200 & AIC  & 1.6678  &  1.8637  &  2.1852  & & 10 (0) & 14 (0) & 16 (0) & & 0  & &  1  \\
    & BIC  & 1.6598  &  1.8637   & 2.1852   & & 9.5 (0) & 13.5 (0) & 16 (0) & & 0  & &  1 \\
    & GAIC & 1.1860  &  1.3045 &   1.4507   & & 0 (0) & 1 (0) & 3 (0)  & & 0.27  & &  0.99  \\
    & GBIC & 1.5115  &  1.7947  &  2.1321   & & 7 (0) & 11 (0) &  15 (0)  & & 0.03  & &  1 \\
    & $\mbox{GBIC}_p$ & 1.1561  &  1.2317 &   1.3614   & & 0 (0) &  1 (0) & 2 (0)  & & 0.43  & &  0.99\\
    & Oracle & 1.1201   & 1.1661 &   1.2167   & &   0 (0)   &  0 (0)   &  0 (0)  & & 1  & &  1 \\ \hline
\end{tabular}
\label{tab5.1.2}
\end{table}

To compare the selected models by the information criteria, we consider the oracle estimate based on the oracle working model $\{1, \cdots, 5\}$ as the benchmark and use both measures of prediction and variable selection. The first performance measure is the prediction error defined as $E (Y - \bx\t \hbbeta)^2$ with $\hbbeta$ an estimate and $(\bx\t, Y)$ an independent observation. The second and third measures are the numbers of false positives and false negatives, respectively. Here, a false positive means  a selected covariate outside the oracle working model and a false negative means a missed covariate in the oracle model. The fourth measure is the probability of recovering exactly the oracle working model, and the fifth measure is the inclusion probability of including all covariates in the oracle working model. The former characterizes the model selection consistency property, while the latter characterizes the sure screening property (Fan and Lv, 2008). In the calculation of the expectation in the prediction error, an independent test sample of size 10,000 was generated for approximating the expectation.

Table 2 summarizes the comparison results of the $25$th, $50$th, and $75$th percentiles of the first three performance measures as well as the selection probability and inclusion probability of the estimated model. The conclusions are similar to those in Section \ref{sec5.1.1}. As the dimensionality $p$ increases, all information criteria tend to select a larger model and the prediction performance deteriorates. GAIC performed better than AIC, and GBIC and $\mbox{GBIC}_p$ performed better than BIC in terms of prediction and variable selection. In particular, the model selected by $\mbox{GBIC}_p$ mimics closely the oracle estimate. We observe that although all methods can have good sure screening property, the probability of model selection consistency generally decreases with dimensionality. It is interesting that GAIC tended to perform better than GBIC when the dimensionality became larger, which can be understood as an effect of decreasing effective sample size due to the noise accumulation associated with high dimensionality.

\subsubsection{High-dimensional linear regression with interaction and moderate effects} \label{sec5.1.3}

\begin{table}
\caption{Percentiles of prediction error, false positives (FP), and false negatives (FN) as well as selection probability (SP) and inclusion probability (IP) of the SICA estimate for model in Section \ref{sec5.1.3} over $100$ simulations for different $p$}
\centering
\begin{tabular}{clccccccccccc}
\hline
$p$ & Criterion & \multicolumn{3}{c}{Prediction error} & & \multicolumn{3}{c}{FP (FN)} & & SP & & IP \\
\cline{3-5} \cline{7-9}
 & & 25th & 50th & 75th & & 25th & 50th & 75th & &   & &   \\
\cline{3-5} \cline{7-9} \cline{11-11} \cline{13-13}
50 & AIC  & 1.3743  &  1.5164  &  1.7514  & & 3 (0)  &  5 (0)   &  7 (0)   & & 0.01  & &  0.94  \\
    & BIC  & 1.2867  &  1.4375  &  1.6424   & & 0 (0)   &  2 (0)  &   4 (0) & & 0.22  & &  0.79   \\
    & GAIC & 1.2427  & 1.3988 &   1.6279   & & 0 (0) &   0 (0) &  1 (1)  & & 0.36  & &  0.61 \\
    & GBIC & 1.2586 &   1.3859   & 1.5947   & & 0 (0) &    1 (0)  &   3 (0)  & & 0.32  & &  0.78 \\
    & $\mbox{GBIC}_p$ & 1.2359  & 1.3793   & 1.6951   & & 0 (0)  &   0 (0)    & 1 (1)  & & 0.38  & &  0.54 \\
    & Oracle & 1.1723  &  1.2481  &  1.3116   & &   0 (0)   &  0 (0)   &  0 (0)  & & 1  & &  1  \\
\cline{3-5} \cline{7-9} \cline{11-11} \cline{13-13}
100 & AIC  & 1.5881  &  1.7814   & 2.0827  & & 6 (0) & 9 (0) & 11 (0)   & & 0  & &  0.85 \\
    & BIC  & 1.5572 &  1.7290 &  2.0359   & & 5 (0) & 8 (0) & 11 (0) & & 0.04  & &  0.80 \\
    & GAIC & 1.2877  &  1.4812  &  1.8517  & & 0 (0) &  0 (1) &  1 (2) & & 0.29  & &  0.42 \\
    & GBIC & 1.3751  &  1.6393   & 1.8202   & & 1 (0) & 4 (0) &  8 (1) & & 0.15  & &  0.74 \\
    & $\mbox{GBIC}_p$ & 1.2877 &   1.4626  &  1.9309   & & 0 (0)  & 0 (1) &  1 (2) & & 0.29  & & 0.42 \\
    & Oracle & 1.1837  &  1.2405  & 1.3155   & &   0 (0)   &  0 (0)   &  0 (0)  & & 1  & &  1 \\
\cline{3-5} \cline{7-9} \cline{11-11} \cline{13-13}
200 & AIC  & 1.7916  &  2.0573  &  2.4881  & & 9.5 (0) & 12 (0) & 15 (0) & & 0  & &  0.78  \\
    & BIC  & 1.7735  &  2.0465   & 2.4881   & & 9 (0) & 12 (0) & 15 (0) & & 0  & &  0.78 \\
    & GAIC & 1.3964  &  1.8119 &   2.2214   & & 0 (0.5) & 0 (2) & 1 (3)  & & 0.21  & &  0.25  \\
    & GBIC &  1.7128  &  2.0465  &  2.4881   & & 9 (0) & 12 (0) &  15 (0)  & & 0.01  & &  0.77 \\
    & $\mbox{GBIC}_p$ & 1.3683  &  1.6990 &   2.1528   & & 0 (0) &  0 (1) & 1 (3)  & & 0.22  & &  0.27\\
    & Oracle & 1.1855   & 1.2367 &   1.3168   & &   0 (0)   &  0 (0)   &  0 (0)  & & 1  & &  1 \\ \hline
\end{tabular}
\label{tab5.1.3}
\end{table}

We now consider the case of high-dimensional linear regression with interaction and moderate effects. We simulated 100 data sets from model (\ref{092}) with the same setting as in Section \ref{sec5.1.2}, except that we set $\bbeta$ as $\bbeta_0 = (1, -1.25, 0.75, -0.95, 1.5, 0.5, -0.5, 0.5, -0.5, 0.5, \cdots, 0)\t$. Compared with the model in Section \ref{sec5.1.2}, the five additional main linear effects are moderately strong. For each data set, we fit the misspecified linear regression model (\ref{053}) without interaction. As argued in Section \ref{sec5.1.2}, the oracle working model is $\supp(\bbeta_0) = \{1, \cdots, 10\}$. We applied the same method to build a sequence of sparse candidate models and used the same performance measures to compare the five model selection criteria as in Section \ref{sec5.1.2}. Table 3 summarizes the comparison results of the $25$th, $50$th, and $75$th percentiles of the first three performance measures as well as the selection probability and inclusion probability of the estimated model. The conclusions are similar to those in Section \ref{sec5.1.2}. We observe that due to the appearance of many moderate effects, all information criteria have less probability of sure screening especially when the dimensionality is large, meaning that some covariates in the oracle working model can be missed.

\subsection{Nonlinear models} \label{sec5.2}
In this section, we examine behaviors of the five information criteria in nonlinear models. They include a polynomial model, a single-index model, and a high-dimensional logistic regression model with interaction. We again observed that the new information
criteria GBIC and GBIC$_p$ outperformed the classical AIC and BIC in all settings.

\subsubsection{Polynomial regression with heteroscedasticity} \label{sec5.2.1}
Since each of the criteria GAIC, GBIC, and $\mbox{GBIC}_p$ includes a second-order term that is related to the variance of the
response, we consider the cubic polynomial regression model with heteroscedastic variance
\begin{equation} \label{094}
y = 1 + 5 x - 2 x^2 + 1.55 x^3 + |x|^{1/2} \varepsilon,
\end{equation}
where $y$ is the response, $x \sim N(0,1)$ is the covariate, and the error $\varepsilon \sim N(0, \sigma^2)$ is independent of $x$. We simulated 100 data sets from this model, each of which contains $n$ i.i.d. observations, with $n = 20$, $40$, $80$, and $\sigma = 0.25$ and $0.5$, respectively. Although the data were generated from model (\ref{094}), we fit the polynomial regression model with constant variance which is incorrectly specified. For each data set, we first fit the polynomial regression model with constant variance of order 1 up to 6 and then selected the final model using the same five model selection criteria. Table 4 summarizes the comparison
results of the frequency of estimated order of polynomial regression model. The conclusions are similar to those in Section \ref{sec5.1.1}. The improved performance of GAIC, GBIC, and $\mbox{GBIC}_p$ over AIC and BIC indicates the usefulness of the second-order terms in capturing the difference between the estimated error variance and the apparent residual variance.

\begin{table}
\caption{Frequency of estimated order for model in Section \ref{sec5.2.1} over $100$ simulations for different $n$ and $\sigma$}
\centering
\begin{tabular}{clccccccccccccccccc}
\hline
$\sigma$ & Criterion & \multicolumn{5}{c}{Model size ($n = 20$)} & & \multicolumn{5}{c}{Model size ($n = 40$)} & & \multicolumn{5}{c}{Model size ($n = 80$)} \\
\cline{3-7} \cline{9-13} \cline{15-19}
 & & 2 & 3 & 4 & 5 & 6 & & 2 & 3 & 4 & 5 & 6 & & 2 & 3 & 4 & 5 & 6 \\
\cline{3-7} \cline{9-13} \cline{15-19}
0.25 & AIC  & 0  &  48  &  17  & 19   & 16 & & 0  &  47  &  23 &   14  &  16 & & 0  &  54  &  17  &  12   & 17 \\
    & BIC  & 0  &  62  &  13 &   14  &  11 & & 0  &  63  &  21  &   9 &    7 & & 0  &  81  &  10 &   4  &   5 \\
    & GAIC & 5   & 78  &  5 &   2  &  4 & & 3  &  88   & 3 &   4  &  2 & & 0   & 97   & 2 &   1 &   0 \\
    & GBIC & 0  &  97   &  3   &  0   &  0 & & 0  &  92   &  5  &  1  &   2 & & 0  &  91  &   5  &   2  &   2 \\
    & $\mbox{GBIC}_p$  & 3   & 94  &   0  &   0 &   0 & & 1  &  98 &    1  &   0  &   0 & & 0  &  99 &    1   &  0  &   0 \\
\cline{3-7} \cline{9-13} \cline{15-19}
0.5   & AIC  & 0  &  38 &   16  &  23 &   23 & & 0  &  40  &  16 &   20   & 24 & & 0  &  47 &   13  &  21  &  19 \\
    & BIC  & 0   & 55 &   15 &   18  &  12 & & 0   & 64 &   14 &   12   & 10 & & 0   & 81   &  9 &   5 &   5 \\
    & GAIC & 12  &  72 &   1  &  1  &  1 & & 5   & 89  &  4 &   0  &  0 & & 1  &  98 &  1 &  0 &   0 \\
    & GBIC & 0  &  97   &  2   &  1  &  0 & & 0  &  88   & 5  &  4  &   3 & & 0 &   93  &   4 &    2  &   1 \\
    & $\mbox{GBIC}_p$  & 8 &   82  &   0 &   0 &   0 & & 3   & 95   & 1   &  0  &  0 & & 0  &  100 &    0  &   0  &  0 \\
    \hline
\end{tabular}
\label{tab5.2.1}
\end{table}

\subsubsection{Nonlinear regression} \label{sec5.2.2}
As another example of misspecified models, we consider a single-index model in which the response depends on a linear combination of covariates through a nonlinear link function. We simulated 100 data sets from the nonlinear regression model
\begin{equation} \label{095}
\by = \bff(\bX \bbeta) + \bveps,
\end{equation}
where $\by$ is an $n$-dimensional response vector, $\bff(\bz) = (f(z_1), \cdots, f(z_n))\t$ with $f(z) = z^3/(0.5 + z^2)$ and $\bz = (z_1, \cdots, z_n)\t$, $\bX$ is an $n \times p$ design matrix, $\bbeta$ is a $p$-dimensional regression coefficient vector, and the error vector $\bveps \sim N(\bzero,
\sigma^2 I_n)$ is independent of $\bX$. The rest of the setting is the same as that in Section \ref{sec5.1.1} except that as in Section \ref{sec5.1.2}, we consider the case that all covariates are independent of each other and thus set $\Sigma_0 = I_p$. Although the data were generated from model (\ref{095}), we fit the linear regression model (\ref{053}). Therefore, the family of working models excludes the true model due to the nonlinearity. As argued in Section \ref{sec5.1.2}, the oracle working model is $\supp(\bbeta_0) = \{1, 2, 3\}$ which indicates the linear regression model (\ref{053}) with the first three covariates. For each data set, we first applied the best subset linear regression to build all possible submodels of $\{1, \cdots, 6\}$ and then selected the final model using the same five model selection criteria. Table 5 summarizes the comparison results of the frequency of estimated model size with $3^\ast$ representing the oracle working model $\{1, 2, 3\}$ for compactness. The conclusions are similar to those in Section \ref{sec5.1.2}.

\begin{table}
\caption{Frequency of estimated model size for model in Section \ref{sec5.2.2} over $100$ simulations for different $n$ and $a$; $3^\ast$ corresponds to the oracle working model $\{1, 2, 3\}$}
\centering
\begin{tabular}{clccccccccccccccccc}
\hline
$\sigma$ & Criterion & \multicolumn{5}{c}{Model size ($n = 20$)} & & \multicolumn{5}{c}{Model size ($n = 40$)} & & \multicolumn{5}{c}{Model size ($n = 80$)} \\
\cline{3-7} \cline{9-13} \cline{15-19}
 & & 2 & $3^\ast$ & 4 & 5 & 6 & & 2 & $3^\ast$ & 4 & 5 & 6 & & 2 & $3^\ast$ & 4 & 5 & 6 \\
\cline{3-7} \cline{9-13} \cline{15-19}
0.25 & AIC  & 0  &  51  &  41   & 8  &  0 & & 0  &  60   & 32 &    7   &  1 & & 0  &  64  &  32   &  3 &   1 \\
    & BIC  & 0  &  71  &  27  &   2  &   0 & & 0  &  82  &  16  &  2   &  0 & & 0  &  90 &   10  &   0   &  0 \\
    & GAIC & 0  &  81  &  16  &  3  &  0 & & 0   & 74  &  22 &    3   & 1 & & 0  &  75   & 22  &   2  &   1 \\
    & GBIC & 0  &  79 &   20  &   1  &   0 & & 0  &  83   & 15 &    2    & 0 & & 0  &  91  &   9  &   0   &  0 \\
    & $\mbox{GBIC}_p$  & 0  &  94 &   6  &   0   &  0 & & 0  &  91  &  9   &  0   &  0 & & 0  &  96  &   4  &   0   &  0 \\
\cline{3-7} \cline{9-13} \cline{15-19}
0.5   & AIC  & 0  &  52 &   36  &  11  &   1 & & 0   & 54  &  33 &   12  &  1 & & 0  &  51   & 40  &   9 &    0 \\
    & BIC  & 0  &  66  &  28   &  5  &   1 & & 0  &  77  &  21  &   2 &    0 & & 0 &   93   &  6   &  1  &   0 \\
    & GAIC & 3  & 73 &   22  &  1  &   0 & & 0 &   67  &  27 &   6 &  0 & & 0   & 63 &   30  &  7 &   0 \\
    & GBIC & 0  &  77  &  20  &   3 &  0 & & 0  &  83  &  17  &   0  &   0 & & 0  &  94  &   6  &   0   &  0 \\
    & $\mbox{GBIC}_p$  & 1  &  90   & 9  &   0  &  0 & & 0 &   90 &   10  &   0  &   0 & & 0 &   99 &    1   &  0 &    0 \\
    \hline
\end{tabular}
\label{tab5.2.2}
\end{table}

\begin{table}
\caption{Percentiles of prediction error, false positives (FP), and false negatives (FN) as well as selection probability (SP) and inclusion probability (IP) of the SICA estimate for model in Section \ref{sec5.2.3} over $100$ simulations}
\centering
\begin{tabular}{lccccccccccc}
\hline
Criterion & \multicolumn{3}{c}{Prediction error} & & \multicolumn{3}{c}{FP (FN)} & & SP & & IP \\
\cline{2-4} \cline{6-8}
 & 25th & 50th & 75th & & 25th & 50th & 75th & &   & &   \\
\cline{2-4} \cline{6-8} \cline{10-10} \cline{12-12}
AIC  & 0.1411  &  0.1512  &  0.1623  & & 15 (0)  &  17 (0)   &  19 (0)   & & 0  & &  0.98  \\
BIC  & 0.1330  &  0.1468  &  0.1591   & & 12.5 (0)   &  15 (0)  &   18 (0) & & 0  & &  0.98   \\
GAIC & 0.1148  & 0.1206 &  0.1284   & & 9 (0) &   11 (0) &  14 (0)  & & 0  & &  0.99 \\
GBIC & 0.1017 &   0.1063   & 0.1155   & & 1 (0) &    2 (0)  &   5.5 (0)  & & 0.23  & &  0.96 \\
$\mbox{GBIC}_p$ & 0.1010  & 0.1048   & 0.1094   & & 0 (0)  &   1 (0)    & 2 (0)  & & 0.39  & & 0.94 \\
Oracle & 0.0923  &  0.0944  &  0.0968   & &   0 (0)   &  0 (0)   &  0 (0)  & & 1  & &  1  \\
\hline
\end{tabular}
\label{tab5.2.3}
\end{table}

\subsubsection{High-dimensional logistic regression with interaction} \label{sec5.2.3}
We now consider the high-dimensional logistic regression with interaction. We simulated 100 data sets from the logistic regression model (\ref{003}) with interaction and an $n$-dimensional parameter vector
\begin{equation} \label{096}
\btheta = \bX \bbeta + 2 \bx_{p + 1},
\end{equation}
where $\bX = (\bx_1, \cdots, \bx_p)$ is an $n \times p$ design matrix, $\bx_{p+1} = \bx_1 \circ \bx_2$ is an interaction term, and the rest is the same as in (\ref{092}). For each data set, the $n$-dimensional response vector $\by$ was sampled from the Bernoulli distribution with success probability vector $(e^{\theta_1}/(1+e^{\theta_1}), \cdots, e^{\theta_n}/(1+e^{\theta_n}))\t$, where $\btheta = (\theta_1, \cdots, \theta_n)\t$ is given in (\ref{096}). As in Section \ref{sec5.1.2}, we consider the case that all covariates are independent of each other and thus set $\Sigma_0 = I_p$. We chose $\bbeta$ as $\bbeta_0 = (2.5, -1.9, 2.8, -2.2, 3, 0, \cdots, 0)\t$ and considered $(n, p) = (200, 1000)$. Although the data were generated from model (\ref{003}) with parameter vector (\ref{096}), we fit the logistic regression model without interaction. This provides another example of misspecified models. As argued in Section \ref{sec5.1.2}, the oracle working model is $\supp(\bbeta_0) = \{1, \cdots, 5\}$ which indicates the logistic regression model (\ref{003}) with the first five covariates. As in Section \ref{sec5.1.2}, for each data set, we first applied the SICA (Lv and Fan, 2009) implemented with the ICA algorithm (Fan and Lv, 2011) to build a sequence of sparse models and then selected the final model using the same five model selection criteria.

We use the same five performance measures as in Section \ref{sec5.1.2}. The prediction error is defined as $E [Y - \exp(\bx\t \hbbeta)/(1+\exp(\bx\t \hbbeta))]^2$ with $\hbbeta$ an estimate and $(\bx\t, Y)$ an independent observation. Table 6 summarizes the comparison results of the $25$th, $50$th, and $75$th percentiles of the first three performance measures as well as the selection probability and inclusion probability of the estimated model. The conclusions are similar to those in Section \ref{sec5.1.2}. GAIC, GBIC, and $\mbox{GBIC}_p$ improved over AIC and BIC in terms of prediction and variable selection.

\section{Discussions} \label{sec6}
We have considered the problem of model selection in misspecified models under two well-known model selection principles: the Bayesian principle and the Kullback-Leibler divergence principle. The novel asymptotic expansions of the two principles in misspecified
generalized linear models lead to the GBIC and GAIC. In particular, a specific form of prior probabilities motivated by the Kullback-Leibler divergence principle leads to the $\mbox{GBIC}_p$ which has a natural decomposition into the negative maximum quasi-log-likelihood and two penalties on model dimensionality and model misspecification, respectively. Numerical studies have demonstrated the advantage of the new methods for model selection in both correctly specified and misspecified models.

When the model is misspecified, the posterior distribution of the parameter may not reflect the true uncertainty of the estimation. In particular, the QMLE $\hbbeta_n$ can have a different asymptotic covariance matrix than the posterior covariance matrix of $\bbeta$. When this difference grows larger as the sample size $n$ increases, its impact on model selection can be rather significant. To mitigate this issue, we have considered an asymptotic framework that allows the prior distribution to vary with sample size, and introduced a local constraint on the prior density preventing it from giving diminishingly low probability to a shrinking region containing the theoretically best model.

Our technical analysis reveals that the covariance contrast matrix, $\bH_n = \bA_n^{-1} \bB_n$, between the covariance structures in the misspecified model and in the true model plays a pivotal role in our model selection framework. Specifically, effective estimates of its trace and log-determinant are needed for implementing the model selection criteria. Through numerical studies, we found that estimation of these two quantities exhibits different behavior. We have considered the simple estimator for $\log |\bH_n|$ and the bootstrap estimator for $\tr\left(\bH_n\right)$; these estimators may not be the best possible ones and the estimation of both quantities $\tr\left(\bH_n\right)$ and $\log |\bH_n|$ needs rigorous theoretical studies, especially in high dimensions. Some other possible ways of estimating $\bH_n$ have been discussed in Section \ref{sec3.4}. The problem of model selection in misspecified models is challenging and our studies can be regarded as the first step toward understanding the effect of explicitly incorporating model misspecification for model selection.

So far we have explored the expression of $\mbox{GBIC}_p$ taking an additive form of the three terms: goodness of fit, model complexity, and model misspecification. Possible extensions of the $\mbox{GBIC}_p$ include introducing other weights and considering non-additive forms. We considered the classical asymptotic setting of fixed dimensionality in our asymptotic expansions for technical simplicity. It would be interesting to extend all the results to high dimensions as well as more general model settings, with the dimensionality diverging with the sample size. It also remains open to characterize the optimality of these model selection criteria in misspecified models. These problems are beyond the scope of the current paper and will be interesting topics for future research.

\appendix

\section{Proofs} \label{secA}
The proofs of Theorem \ref{thm4} and Propositions \ref{prop1} and \ref{prop2} are omitted here and are given in the technical report Lv and Liu (2010).

\subsection{Proof of Theorem \ref{thm1}} \label{secA.1}
We condition on the event $\widetilde{Q}_n = \{\hbbeta_n \in N_n(\delta_n)\}$, where $N_n(\delta_n) = \{\bbeta \in \mathbb{R}^d: \|\bB_n^{1/2} (\bbeta - \bbeta_{n, 0})\| \leq \delta_n\}$. By the proof of consistency of the QMLE $\hbbeta_n$ in Theorem \ref{thm5}, we have
\[ P(\widetilde{Q}_n) \geq 1 - \frac{4 d}{c^2 \delta_n^2} \rightarrow 1 \]
since $\delta_n \rightarrow \infty$ by Condition \ref{cond5}. It is easy to see that $\ell^*_n(\by, \bbeta)$ defined in (\ref{084}) is a smooth concave function on $\mathbb{R}^d$ with its maximum 0 attained at $\bbeta = \hbbeta_n$,
\[ \partial^2 \ell^*_n(\by, \bbeta)/\partial \bbeta^2 = -\bA_n(\bbeta), \]
and $U_n(\bbeta) = \exp\left[n^{-1} \ell^*_n(\by, \bbeta)\right]$ takes values between 0 and 1. Thus by Taylor's theorem, expanding $\ell^*_n(\by, \cdot)$ around $\hbbeta_n$ gives for any $\bbeta$ in the neighborhood $\widetilde{N}_n(\delta_n) = \{\bbeta \in \mathbb{R}^d: \|\bB_n^{1/2} (\bbeta - \hbbeta_n)\| \leq \delta_n\}$,
\begin{equation} \label{026}
\ell^*_n(\by, \bbeta) = \frac{1}{2} (\bbeta - \hbbeta_n)\t \left[\partial^2 \ell^*_n(\by, \bbeta_*)/\partial \bbeta^2\right] (\bbeta - \hbbeta_n) = -\frac{n}{2} \bdelta\t \bV_n(\bbeta_*) \bdelta,
\end{equation}
where $\bbeta_*$ lies on the line segment joining $\bbeta$ and $\hbbeta_n$, $\bdelta = n^{-1/2} \bB_n^{1/2} (\bbeta - \hbbeta_n)$, and $\bV_n(\bbeta) = \bB_n^{-1/2} \bA_n(\bbeta) \bB_n^{-1/2}$. Clearly, the neighborhood $\widetilde{N}_n(\delta_n)$ defined above is a convex set containing $\hbbeta_n$ and thus $\bbeta_* \in \widetilde{N}_n(\delta_n) \subset N_n(2 \delta_n)$. Let
\[ \rho_n(\delta_n) = \max_{\bbeta \in N_n(2 \delta_n)} \max\left\{|\lambda_{\min}\left(\bV_n(\bbeta) - \bV_n\right)|, |\lambda_{\max}\left(\bV_n(\bbeta) - \bV_n\right)|\right\} \]
with $\bV_n = \bV_n(\bbeta_{n, 0})$. By Conditions \ref{cond5} and \ref{cond2}, $\rho_n(\delta_n) \leq \lambda_{\min}(\bV_n)/2$ for sufficiently large $n$. Thus it follows from (\ref{026}) that
\begin{equation} \label{025}
q_1(\bbeta) 1_{\widetilde{N}_n(\delta_n)}(\bbeta) \leq -n^{-1} \ell_n^*(\by, \bbeta) 1_{\widetilde{N}_n(\delta_n)}(\bbeta) \leq q_2(\bbeta) 1_{\widetilde{N}_n(\delta_n)}(\bbeta),
\end{equation}
where $q_1(\bbeta) = \frac{1}{2} \bdelta\t \left[\bV_n - \rho_n(\delta_n) I_d\right] \bdelta$ and $q_2(\bbeta) = \frac{1}{2} \bdelta\t \left[\bV_n + \rho_n(\delta_n) I_d\right] \bdelta$.
This entails
\begin{equation} \label{027}
E_{\mu_\mathfrak{M}} \left(e^{-n q_2} 1_{\widetilde{N}_n(\delta_n)}\right) \leq E_{\mu_\mathfrak{M}} \left[U_n(\bbeta)^n 1_{\widetilde{N}_n(\delta_n)}\right] \leq E_{\mu_\mathfrak{M}} \left(e^{-n q_1} 1_{\widetilde{N}_n(\delta_n)}\right).
\end{equation}
We will see later that inequality (\ref{027}) is the key step in deriving the asymptotic expansion of $\log E_{\mu_\mathfrak{M}} [U_n(\bbeta)^n]$ in (\ref{023}).

We list some auxiliary results in the following two lemmas, whose proofs are similar to those given in the technical report Lv and Liu (2010).

\begin{lemma} \label{lem1}
Under Condition \ref{cond5} in Section \ref{sec4.1}, we have
\begin{equation}
\label{030}
c_1 \int e^{-n q_j} 1_{\widetilde{N}_n(\delta_n)} d\mu_0 \leq E_{\mu_\mathfrak{M}} \left(e^{-n q_j} 1_{\widetilde{N}_n(\delta_n)}\right) \leq c_2 \int e^{-n q_j} 1_{\widetilde{N}_n(\delta_n)} d\mu_0, \quad j = 1, 2.
\end{equation}
\end{lemma}

\begin{lemma} \label{lem2}
It holds that
\begin{align}
\label{028}
E_{\mu_\mathfrak{M}} \left[U_n(\bbeta)^n 1_{\widetilde{N}_n^c(\delta_n)}\right] & \leq \exp\left(-\kappa_n \delta_n^2\right), \\
\label{031}
\int e^{-n q_1} d\mu_0 & = \left(\frac{2 \pi}{n}\right)^{d/2} |\bV_n - \rho_n(\delta_n) I_d|^{-1/2}, \\
\label{032}
\int e^{-n q_2} d\mu_0 & = \left(\frac{2 \pi}{n}\right)^{d/2} |\bV_n + \rho_n(\delta_n) I_d|^{-1/2}, \\
\label{029}
\int e^{-n q_j} 1_{\widetilde{N}_n^c(\delta_n)} d\mu_0 & \leq  \left(\frac{2 \pi}{n \kappa_n}\right)^{d/2} \exp\left[-\frac{1}{2} \kappa_n \delta_n^2 + \frac{d}{2} + \frac{d}{2} \log \left(\kappa_n \delta_n^2 d^{-1}\right)\right], j = 1, 2,
\end{align}
where $\kappa_n = \lambda_{\min}(\bV_n)/2$ and it is assumed that $\kappa_n  \delta_n^2 > d$ in (\ref{029}).
\end{lemma}

Now we are ready to obtain the asymptotic expansion of $\log E_{\mu_\mathfrak{M}}
[U_n(\bbeta)^n]$ in (\ref{023}). By Conditions \ref{cond5} and \ref{cond2}, $\lambda_{\min}^{-1}(\bV_n) \delta_n^{-2} = o(\log^{-1} n)$, which along with (\ref{028}) shows that the term $E_{\mu_\mathfrak{M}} [U_n(\bbeta)^n 1_{\widetilde{N}_n^c(\delta_n)}]$ converges to zero at a rate faster than any polynomial rate in $n$. Similarly, in view of (\ref{029}), we can show that the two terms $\int e^{-n q_j} 1_{\widetilde{N}_n^c(\delta_n)} d\mu_0$ with $q = 1, 2$ converge to zero at a rate faster than any polynomial rate in $n$. It follows from $\rho_n(\delta_n) = o(\lambda_{\min}(\bV_n))$ that
\begin{align*}
|\bV_n \pm \rho_n(\delta_n) I_d|^{-1/2} & = |\bV_n|^{-1/2} |I_d \pm \rho_n(\delta_n) \bV_n^{-1}|^{-1/2} \\
& = |\bV_n|^{-1/2}
\left\{1 + O\left[\rho_n(\delta_n) \tr(\bV_n^{-1})\right]\right\} \\
& = |\bV_n|^{-1/2}
\left\{1 + O\left[\rho_n(\delta_n) d \lambda_{\min}^{-1}(\bV_n)\right]\right\} =  |\bV_n|^{-1/2}
[1 +o(1)].
\end{align*}
Since $\lambda_{\max}(\bV_n)$ is assumed to be of a polynomial order of $n$ in Condition \ref{cond5}, combining these results with (\ref{027}) and Lemmas \ref{lem1} and \ref{lem2} yields
\begin{align*}
\log E_{\mu_\mathfrak{M}} [U_n(\bdelta)^n] & = \log \left\{ \left(\frac{2
\pi}{n}\right)^{d/2} |\bV_n|^{-1/2} [1 +o(1)]\right\} + \log c_n \\
& = -\frac{\log n}{2} d + \frac{1}{2} \log |\bA_n^{-1}
\bB_n| + \frac{\log (2 \pi)}{2} d +\log c_n + o(1),
\end{align*}
where $c_n \in [c_1, c_2]$. This together with (\ref{023}) proves the conclusion on the event $\widetilde{Q}_n$.

\subsection{Proof of Theorem \ref{thm2}} \label{secA.2}
By (\ref{006}) and the first part of the proof of Theorem \ref{thm3}, we have
\[ \log \alpha_\mathfrak{M} =  -\frac{1}{2} \tr(\bH_n) + \log C + o(1), \]
which along with Theorem \ref{thm1} leads to the desired asymptotic expansion of the Bayes factor.

\subsection{Proof of Theorem \ref{thm3}} \label{secA.3}
It suffices to prove
\begin{equation} \label{047}
E \eta_n(\hbbeta_n) = \eta_n(\bbeta_{n, 0}) - \frac{1}{2} \tr\left(\bA_n^{-1} \bB_n\right) + o(1)
\end{equation}
and
\begin{equation} \label{048}
\eta_n(\bbeta_{n, 0}) = E \left[\ell_n(\by, \hbbeta_n)\right] - \frac{1}{2} \tr\left(\bA_n^{-1} \bB_n\right) + o(1).
\end{equation}
We will prove (\ref{047}) and (\ref{048}) in separate steps. Define
\[ \widetilde{\ell}_n(\by, \bbeta) = \ell_n(\by, \hbbeta_n + n^{1/2} \bC_n^{-1} \bbeta) \quad \text{ and } \quad \widetilde{\eta}_n(\bbeta) = \eta_n(\bbeta_{n, 0} + n^{1/2} \bC_n^{-1} \bbeta), \]
where $\bC_n = \bB_n^{-1/2} \bA_n$ and $\bbeta = (\beta_1, \cdots, \beta_d)\t$.

1) We first prove (\ref{047}). The key idea is to do a second-order Taylor expansion of $\widetilde{\eta}_n(\bbeta)$ around $\bzero$ and retain the Lagrange remainder term. By Theorem \ref{thm4}, $\eta_n$ attains its maximum at $\bbeta_{n, 0}$. Thus we derive
\begin{align*}
\eta_n(\hbbeta_n) & = \widetilde{\eta}_n\left[n^{-1/2} \bC_n (\hbbeta_n - \bbeta_{n, 0})\right] = \widetilde{\eta}_n(\bzero) - \frac{n}{2} \bv\t \left(\bC_n^{-1} \bA_n \bC_n^{-1}\right) \bv \\
& \quad - \frac{1}{6} \sum_{j, k, l} \left[\partial^3 \widetilde{\eta}_n(\bbeta_*)/\partial \beta_j \partial \beta_k \partial \beta_l\right] v_j v_k v_l \\
& = \eta_n(\bbeta_{n, 0}) - \frac{n}{2} \bv\t \left(\bB_n^{1/2} \bA_n^{-1} \bB_n^{1/2}\right) \bv - \frac{1}{6} \sum_{j, k, l} \left[\partial^3 \widetilde{\eta}_n(\bbeta_*)/\partial \beta_j \partial \beta_k \partial \beta_l\right] v_j v_k v_l,
\end{align*}
where $\bv = (v_1, \cdots, v_d)\t = n^{-1/2} \bC_n (\hbbeta_n - \bbeta_{n, 0})$ and $\bbeta_*$ lies on the line segment joining $\bv$ and $\bzero$. By Condition \ref{cond6},
\[ |\partial^3 \widetilde{\eta}_n(\bbeta_*)/\partial \beta_j \partial \beta_k \partial \beta_l| = o(n^{3/2}). \]
Since $\bC_n (\hbbeta_n - \bbeta_{n, 0}) \toD N(\bzero, I_d)$ by Theorem \ref{thm6} and $E \|\bC_n (\hbbeta_n - \bbeta_{n, 0})\|^{3 + \alpha} = O(1)$ for some $\alpha > 0$ by Condition \ref{cond6}, Theorem 6.2 in DasGupta (2008) applies to show that for any $j, k, l$,
\[ E (n v_j v_k) \rightarrow \delta_{jk} \quad \text{ and } \quad E(n^{3/2} |v_j v_k v_l|) = O(1), \]
where $\delta_{jk}$ denotes the Kronecker delta. By Condition \ref{cond2}, we have $\lambda_{\min}(\bV_n) \geq c$ for some positive constant $c$, which entails that
\begin{align*}
\tr\left(\bA_n^{-1} \bB_n\right) & = \tr\left(\bB_n^{1/2} \bA_n^{-1} \bB_n^{1/2}\right) \leq d \lambda_{\max}(\bB_n^{1/2} \bA_n^{-1} \bB_n^{1/2}) \\
& =  d \lambda_{\min}^{-1}(\bV_n) = O(1).
\end{align*}
Combining these results yields
\begin{align*}
E \eta_n(\hbbeta_n) & = \eta_n(\bbeta_{n, 0}) - \left[\frac{1}{2} \tr \left(\bB_n^{1/2} \bA_n^{-1} \bB_n^{1/2}\right) + o(1)\right] + \frac{1}{6} \sum_{j, k, l} o(n^{3/2}) O(n^{-3/2}) \\
& = \eta_n(\bbeta_{n, 0}) - \frac{1}{2} \tr\left(\bA_n^{-1} \bB_n\right) + o(1).
\end{align*}

2) We then prove (\ref{048}). The key idea is to represent $\ell_n(\by, \bbeta_{n, 0})$ using a second-order Taylor expansion of $\ell_n(\by, \cdot)$ around $\hbbeta_n$ and retaining the Lagrange remainder term. It is easy to see that the difference between $\eta_n(\bbeta)$ and $\ell_n(\by, \bbeta)$ is linear in $\bbeta$, which implies that each pair of their second or higher order partial derivatives agree. Since $\ell_n(\by, \cdot)$ attains its maximum at $\hbbeta_n$, we derive
\begin{align*}
\ell_n(\by, \bbeta_{n, 0}) & = \widetilde{\ell}_n(\by, \bv) = \widetilde{\ell}_n(\bzero) - \frac{1}{2} \bv\t \left[\partial^2 \widetilde{\ell}_n(\by, \bzero)/\partial \bbeta^2\right] \bv \\
& \quad - \frac{1}{6} \sum_{j, k, l} \left[\partial^3 \widetilde{\ell}_n(\by, \bbeta_*)/\partial \beta_j \partial \beta_k \partial \beta_l\right] v_j v_k v_l \\
& = \ell_n(\by, \hbbeta_n) - \frac{1}{2} \bv\t \left[\partial^2 \widetilde{\eta}_n(-\bv)/\partial \bbeta^2\right] \bv - \frac{1}{6} \sum_{j, k, l} \left[\partial^3 \widetilde{\eta}_n(\bbeta_{**})/\partial \beta_j \partial \beta_k \partial \beta_l\right] v_j v_k v_l,
\end{align*}
where $\bv = (v_1, \cdots, v_d)\t = n^{-1/2} \bC_n (\bbeta_{n, 0} - \hbbeta_n)$, $\bbeta_*$ lies on the line segment joining $\bv$ and $\bzero$, and $\bbeta_{**} = \bbeta_* - \bv$.

It follows from $\sup_{\bbeta} \max_{j, k, l} |\partial^3 \widetilde{\eta}_n(\bbeta)/\partial \beta_j \partial \beta_k \partial \beta_l| = o(n^{3/2})$ in Condition \ref{cond6} that
\[ |\partial^3 \widetilde{\eta}_n(\bbeta_{**})/\partial \beta_j \partial \beta_k \partial \beta_l| = o(n^{3/2}) \]
and $\partial^2 \widetilde{\eta}_n(\bbeta)/\partial \bbeta^2$ is Lipschitz with Lipschitz constant $o(n^{3/2})$. Thus
\begin{align*}
\bv\t \left[\partial^2 \widetilde{\eta}_n(-\bv)/\partial \bbeta^2\right] \bv & = \bv\t \left[\partial^2 \widetilde{\eta}_n(\bzero)/\partial \bbeta^2\right] \bv + o(n^{3/2}) \|\bv\|^3 \\
& = n \bv\t \left(\bC_n^{-1} \bA_n \bC_n^{-1}\right) \bv + o(n^{3/2}) \|\bv\|^3 \\
& = n \bv\t \left(\bB_n^{1/2} \bA_n^{-1} \bB_n^{1/2}\right) \bv + o(n^{3/2}) \|\bv\|^3.
\end{align*}
Since $\bC_n (\hbbeta_n - \bbeta_{n, 0}) \toD N(\bzero, I_d)$ by Theorem \ref{thm6} and $E \|\bC_n (\hbbeta_n - \bbeta_{n, 0})\|^{3 + \alpha} = O(1)$ for some $\alpha > 0$ by Condition \ref{cond6}, Theorem 6.2 in DasGupta (2008) applies to show that for any $j$ and $k$,
\[ E (n v_j v_k) \rightarrow \delta_{jk} \quad \text{ and } \quad E(n^{3/2} \|\bv\|^3) = O(1). \] These along with $\tr\left(\bA_n^{-1} \bB_n\right) = O(1)$ yield
\begin{align*}
\eta_n(\bbeta_{n, 0}) & = E \ell_n(\by, \bbeta_{n, 0}) = E \ell_n(\by, \hbbeta_n) - \left[\frac{1}{2} \tr \left(\bB_n^{1/2} \bA_n^{-1} \bB_n^{1/2}\right) + o(1)\right] + o(n^{3/2}) O(n^{-3/2}) \\
& = E \ell_n(\by, \hbbeta_n) - \frac{1}{2} \tr\left(\bA_n^{-1} \bB_n\right) + o(1),
\end{align*}
which completes the proof.

\subsection{Proof of Theorem \ref{thm7}} \label{secA.6}
Since the model is correctly specified, that is, $G_n = F_n(\cdot, \bbeta_{n, 0})$, we have the identity $\bA_n = \bB_n$, yielding $\bH_n = \bA_n^{-1} \bB_n = I_d$. Note that by assumption, the smallest and largest eigenvalues of $n^{-1} \bB_n$ are bounded away from 0 and $\infty$. To prove the consistency of the estimator $\widehat{\bH}_n = \widehat{\bA}_n^{-1} \widehat{\bB}_n$, it suffices to show that
\[ n^{-1} \widehat{\bA}_n = n^{-1} \bB_n + o_P(1) \quad \text{ and } \quad n^{-1} \widehat{\bB}_n = n^{-1} \bB_n + o_P(1). \]
By the proof of Theorem \ref{thm5}, we have $\bB_n^{1/2} (\hbbeta_n - \bbeta_{n, 0}) = O_P(1)$, which along with the assumption that  $\lambda_{\min}(n^{-1} \bB_n)$ is bounded away from 0 gives
\[ \hbbeta_n = \bbeta_{n, 0} + O_P(n^{-1/2}). \]

By Condition \ref{cond7}, $n^{-1} \bA_n(\bbeta)$ is Lipschitz in an asymptotically shrinking neighborhood of $\bbeta_{n, 0}$ given by $N_n(\delta) = \{\bbeta \in \mathbb{R}^d: \|\bB_n^{1/2} (\bbeta - \bbeta_{n, 0})\| \leq \delta\}$. It follows from the proof of Theorem \ref{thm5} that
\[ P\left\{\hbbeta_n \in N_n(\delta)\right\} \rightarrow 1 \quad \text{ as } n \rightarrow \infty. \]
Combining this result with the Lipschitz continuity of $n^{-1} \bA_n(\bbeta)$ and the consistency of $\hbbeta_n$ yields
\[ n^{-1} \widehat{\bA}_n = n^{-1} \bB_n + o_P(1). \]

It remains to prove $n^{-1} \widehat{\bB}_n = n^{-1} \bB_n + o_P(1)$. To analyze $n^{-1} \widehat{\bB}_n$, we decompose it into three terms
\[
n^{-1} \widehat{\bB}_n = n^{-1} \bX\t \diag\left\{\left[\by - \bmu(\bX
  \hbbeta_n)\right] \circ \left[\by - \bmu(\bX
  \hbbeta_n)\right]\right\} \bX \equiv \bG_1 + \bG_2 + \bG_3,
\]
with
\begin{align*}
\bG_1 & = n^{-1} \bX\t \diag\{(\by - E \by) \circ (\by - E \by)\} \bX, \\
\bG_2 & = -2 n^{-1} \bX\t \diag\{(\by - E \by) \circ [\bmu(\bX \hbbeta_n) - \bmu(\bX \bbeta_{n, 0})]\} \bX, \\
\bG_3 & = n^{-1} \bX\t \diag\{[\bmu(\bX \hbbeta_n) - \bmu(\bX \bbeta_{n, 0})] \circ [\bmu(\bX \hbbeta_n) - \bmu(\bX \bbeta_{n, 0})]\} \bX,
\end{align*}
since $E \by = \bmu(\bX \bbeta_{n, 0})$ by (\ref{073}). Since $\bG_1 = n^{-1} \sum_{i = 1}^n \{\bx_i \bx_i\t (Y_i - E Y_i)^2\}$ is the mean of $n$ independent random variables and has asymptotically vanishing 
second central moment by Condition \ref{cond7}, it follows from the law of large numbers that
\[ \bG_1 = E \bG_1 + o_P(1) = n^{-1} \bB_n + o_P(1). \]
By Condition \ref{cond7}, $n^{-1} \bX\t \diag\{[\bmu(\bX \bbeta) - \bmu(\bX \bbeta_{n, 0})] \circ [\bmu(\bX \bbeta) - \bmu(\bX \bbeta_{n, 0})]\} \bX$ is Lipschitz in the  neighborhood $N_n(\delta)$. Thus similarly as for $n^{-1} \widehat{\bA}_n$, we can show that
\[ \bG_3 = o_P(1). \]

We finally prove $\bG_2 = o_P(1)$. By Condition \ref{cond7}, $Y_i - E Y_i$ are uniformly bounded by $K_n = o(n^{1/2})$ with asymptotic probability one, that is,
\[ P(\max\nolimits_{i = 1}^n |Y_i - E Y_i| \leq K_n) \rightarrow 1 \]
as $n \rightarrow \infty$. Note that on the above event, the largest absolute eigenvalue of the matrix $\bG_2$ is bounded by the largest eigenvalue of
\begin{align*}
2 & n^{-1} \bX\t \diag\{|\by - E \by| \circ |\bmu(\bX \hbbeta_n) - \bmu(\bX \bbeta_{n, 0})|\} \bX \\
& \leq 2 K_n n^{-1} \bX\t \diag\{|\bmu(\bX \hbbeta_n) - \bmu(\bX \bbeta_{n, 0})|\} \bX,
\end{align*}
where the absolute value is understood componentwise and the inequality means the order for symmetric positive semidefinite matrices. By Condition \ref{cond7}, $n^{-1} \bX\t \diag\{|\bmu(\bX \bbeta) - \bmu(\bX \bbeta_{n, 0})|\}$ is Lipschitz in the neighborhood $N_n(\delta)$. Thus by $K_n = o(n^{1/2})$ and $\hbbeta_n = \bbeta_{n, 0} + O_P(n^{-1/2})$, we have $G_2 = o_P(1)$. Combining the above results concludes the proof.

\subsection{Proof of Theorem \ref{thm5}} \label{secA.4}
It is easy to see that $N_n(\delta)$ is a convex set by its definition. Denote by
\[ \partial N_n(\delta) = \left\{\bbeta \in \mathbb{R}^d: \|\bB_n^{1/2} (\bbeta - \bbeta_{n, 0})\| = \delta\right\} \] the boundary of the closed set $N_n(\delta)$. We define an event
\[ Q_n = \left\{\ell_n(\by, \bbeta_{n, 0}) > \max_{\bbeta \in \partial N_n(\delta)} \ell_n(\by, \bbeta)\right\}. \]
We observe that if the event $Q_n$ occurs, the continuous function $\ell_n(\by, \cdot)$ has a local maximum in the interior of $N_n(\delta)$. Since the first part of Condition \ref{cond1} implies that $\ell_n(\by, \cdot)$ is strictly concave, this maximum must be located at $\hbbeta_n$. This shows that
\[ Q_n \subset \{\hbbeta_n \in N_n(\delta)\}. \]

Next we construct a lower bound on $P(Q_n)$. By Taylor's theorem, we have for any $\bbeta$,
\[ \ell_n(\by, \bbeta) - \ell_n(\by, \bbeta_{n, 0}) = (\bbeta - \bbeta_{n, 0})\t \bPsi_n(\bbeta_{n, 0}) - \frac{1}{2} (\bbeta - \bbeta_{n, 0})\t \bA_n(\bbeta_*) (\bbeta - \bbeta_{n, 0}), \]
where $\bbeta_*$ lies on the line segment joining $\bbeta$ and $\bbeta_{n, 0}$. We make a transformation of the variable by letting
\[ \bu = \delta^{-1} \bB_n^{1/2} (\bbeta - \bbeta_{n, 0}). \]
Then it follows that
\[ \ell_n(\by, \bbeta) - \ell_n(\by, \bbeta_{n, 0}) = \delta \bu\t \bB_n^{-1/2} \bPsi_n(\bbeta_{n, 0}) - \delta^2 \bu\t \bV_n(\bbeta_*) \bu/2. \]
Note that $\bbeta \in \partial N_n(\delta)$ if and only if $\|\bu\| = 1$, and that $\bbeta \in \partial N_n(\delta)$ implies $\bbeta_* \in N_n(\delta)$ by the convexity of $N_n(\delta)$. Clearly
\[ \max_{\|\bu\| = 1} \bu\t \bB_n^{-1/2} \bPsi_n(\bbeta_{n, 0}) = \|\bB_n^{-1/2} \bPsi_n(\bbeta_{n, 0})\| \]
and by the second part of Condition \ref{cond2}, for $n$ sufficiently large we have
\[ \min_{\|\bu\| = 1} \bu\t \bV_n(\bbeta_*) \bu \geq \min_{\bbeta \in N_n(\delta)} \lambda_{\min} \left\{\bV_n(\bbeta)\right\} \geq c. \]
Thus we have
\[ \max_{\bbeta \in \partial N_n(\delta)} \ell_n(\by, \bbeta) - \ell_n(\by, \bbeta_{n, 0}) \leq \delta (\|\bB_n^{-1/2} \bPsi_n(\bbeta_{n, 0})\| - c \delta/2), \]
which along with Markov's inequality entails that
\[ P\left(Q_n\right) \geq P\left(\|\bB_n^{-1/2} \bPsi_n\left(\bbeta_{n, 0}\right)\|^2 < c^2 \delta^2/4\right) \geq 1 - \frac{E \|\bB_n^{-1/2} \bPsi_n(\bbeta_{n, 0})\|^2}{c^2 \delta^2/4}. \]

Observe that
\begin{align*}
E \|\bB_n^{-1/2} \bPsi_n(\bbeta_{n, 0})\|^2 & = E \tr\left[\bPsi_n(\bbeta_{n, 0})\t \bB_n^{-1} \bPsi_n(\bbeta_{n, 0})\right] = E \tr\left[\bPsi_n(\bbeta_{n, 0}) \bPsi_n(\bbeta_{n, 0})\t \bB_n^{-1}\right] \\
& = \tr\left\{E\left[\bPsi_n(\bbeta_{n, 0}) \bPsi_n(\bbeta_{n, 0})\t\right] \bB_n^{-1}\right\} = \tr(I_d) = d.
\end{align*}
For any given $\eta \in (0, 1)$, letting $\delta = (2 d^{1/2})/(c \eta^{1/2})$ thus makes
\[ P(Q_n) \geq 1 - \frac{4 d}{c^2 \delta^2} = 1 - \eta. \]
This together with $Q_n \subset \{\hbbeta_n \in N_n(\delta)\}$ and the first part of Condition \ref{cond2} proves $\hbbeta_n - \bbeta_{n, 0} = o_P(1)$, which completes the proof.

\subsection{Proof of Theorem \ref{thm6}} \label{secA.5}
We condition on the event $Q_n = \{\ell_n(\by, \bbeta_{n, 0}) > \max_{\bbeta \in \partial N_n(\delta)} \ell_n(\by, \bbeta)\}$ defined in the proof of Theorem \ref{thm5}, which has been shown to entail that $\hbbeta_n \in N_n(\delta)$. By the mean value theorem applied componentwise and (\ref{011}), we obtain
\begin{align*}
\bzero & = \bPsi_n(\hbbeta_n) = \bPsi_n(\bbeta_{n, 0}) - \widetilde{\bA}_n(\bbeta_1, \cdots, \bbeta_d) (\hbbeta_n - \bbeta_{n, 0}) \\
& = \bX\t (\by - E \by) - \widetilde{\bA}_n(\bbeta_1, \cdots, \bbeta_d) (\hbbeta_n - \bbeta_{n, 0}),
\end{align*}
where each of $\bbeta_1, \cdots, \bbeta_d$ lies on the line segment joining $\hbbeta_n$ and $\bbeta_{n, 0}$. Let $\bC_n = \bB_n^{-1/2} \bA_n$. Then we have
$\bC_n (\hbbeta_n - \bbeta_{n, 0}) = \bu_n + \bw_n$,
where $\bu_n = -\bB_n^{-1/2} \bX\t (\by - E \by)$ and
\[ \bw_n = -\left[\widetilde{\bV}_n(\bbeta_1, \cdots, \bbeta_d) - \bV_n\right] \left[\bB_n^{1/2} (\hbbeta_n - \bbeta_{n, 0})\right]. \]
By Slutsky's lemma and the proof of Theorem \ref{thm5}, we see that to show
$\bC_n (\hbbeta_n - \bbeta_{n, 0}) \toD N(\bzero, I_d)$,
it suffices to prove $\bu_n \toD N(\bzero, I_d)$ and $\bw_n = o_P(1)$.

The result $\bw_n = o_P(1)$ follows from $Q_n \subset \{\hbbeta_n \in N_n(\delta)\}$ and Condition \ref{cond3}. Finally we show that $\bu_n \toD N(\bzero, I_d)$. Fix an arbitrary unit vector $\ba \in \mathbb{R}^d$. We consider the asymptotic distribution of the linear combination
$v_n = \ba\t \bu_n = -\ba\t \bB_n^{-1/2} \bX\t (\by - E \by) = \sum_{i = 1}^n z_i$,
where $z_i = -\ba\t \bB_n^{-1/2} \bx_i (y_i - E y_i)$, $i = 1, \cdots, n$. Clearly $z_i$'s are independent and have mean 0, and
$\sum_{i = 1}^n \var(z_i) = \ba\t \bB_n^{-1/2} \bB_n \bB_n^{-1/2} \ba = 1$.
By Condition \ref{cond4}, $E |y_i - E y_i|^3 \leq M$ for some positive constant $M$. Thus we derive
\begin{align*}
\sum_{i = 1}^n E |z_i|^3 & = \sum_{i = 1}^n |\ba\t \bB_n^{-1/2} \bx_i|^3 E |y_i - E y_i|^3 \leq M \sum_{i=1}^n |\ba\t \bB_n^{-1/2} \bx_i |^3 \\
& \leq M \sum_{i=1}^n \|\ba\|^3\|\bB_n^{-1/2} \bx_i\|^3 = M \sum_{i=1}^n (\bx_i\t \bB_n^{-1} \bx_i)^{3/2}\rightarrow 0,
\end{align*}
where we used the Cauchy-Schwarz inequality and Condition \ref{cond4}. Applying Lyapunov's theorem yields
$\ba\t \bu_n = \sum_{i = 1}^n z_i \toD N(0, 1)$.
Since this asymptotic normality holds for any unit vector $\ba \in \mathbb{R}^d$, we conclude that $\bu_n \toD N(\bzero, I_d)$, which concludes the proof.

\section{Three commonly used GLMs} \label{secB}
In this section, we provide the formulas for matrices $\widehat{\bA}_n$ and $\widehat{\bB}_n$ in (\ref{087})--(\ref{088}), when the working model $F_n$ is chosen to be one of the three commonly used GLMs: the linear regression model, logistic regression model, and Poisson regression model.

\subsection{Linear regression} \label{secB.1}
The linear regression model is given by
\begin{equation} \label{040}
\by = \bX \bgamma + \bveps,
\end{equation}
where $\by$ is an $n$-dimensional response vector, $\bX$ is an $n \times d$ design matrix, $\bgamma$ is a $d$-dimensional regression coefficient vector, and $\bveps \sim N(\bzero, \sigma^2 I_n)$ is an $n$-dimensional error vector. Let $\bbeta = \bgamma/\sigma^2$ and $b(\theta) =
\sigma^2\theta^2/2, \ \theta \in \mathbb{R}$. Then the quasi-log-likelihood of the sample in (\ref{002}) becomes
\begin{align*}
\ell_n(\by, \bbeta) & = \by\t \bX \bbeta - \bone\t \bb(\bX
  \bbeta) - \frac{\|\by\|_2^2}{2 \sigma^2} - \frac{n}{2} \log
\sigma^2 - \frac{n \log 2 \pi}{2} \\
\label{041}
& = -\frac{\|\by - \bX \bgamma\|_2^2}{2 \sigma^2} - \frac{n}{2} \log
\sigma^2 - \frac{n \log 2 \pi}{2}.
\end{align*}
Maximizing $\ell_n(\by, \bbeta)$ with respect to $\bgamma$ yields the least-squares estimator $\hbgamma = (\bX\t \bX)^{-1} \bX \by$.
We estimate the error variance $\sigma^2$ using the residual sum of squares (RSS),
$\widehat{\sigma}^2 = \|\by - \bX \hbgamma\|_2^2/(n - d)$.
Then we obtain an estimate $\hbbeta_n = \hbgamma/\widehat{\sigma}^2$. Since $b'(\theta) = \sigma^2 \theta$ and $b''(\theta) = \sigma^2$, we have
\begin{equation} \label{042}
\widehat{\bA}_n = \widehat{\sigma}^2 \bX\t \bX \quad \text{ and } \quad
\widehat{\bB}_n =  \bX\t \diag\left\{\left[\by - \bX \hbgamma\right] \circ \left[\by - \bX \hbgamma\right]\right\} \bX,
\end{equation}
where $\circ$ denotes the Hadamard (componentwise) product.

An interesting special case is when the true model is indeed linear, but may involve a different set of covariates. Then $\bB_n = \tau^2 \bX^{T}\bX$, where $\tau^2$ is the true variance of the response variable $Y$ had we known the true model precisely. The $\mbox{GBIC}_p$ just involves an additional term penalizing the inflation of the error  variance due to model misspecification.

\subsection{Logistic regression} \label{secB.2}
In the logistic regression model, $b(\theta) = \log(1 + e^\theta)$, $\theta \in \mathbb{R}$ and the quasi-log-likelihood of the sample is
$\ell_n(\by, \bbeta) = \by\t \bX \bbeta - \bone\t \bb(\bX \bbeta)$ with the last constant term in (\ref{002}) ignored, which attains its maximum at $\hbbeta_n$. Since
$b'(\theta) = e^\theta/(1 + e^\theta)$ and $b''(\theta) =
e^\theta/(1 + e^\theta)^2$, matrices $\widehat{\bA}_n$ and $\widehat{\bB}_n$ are given in (\ref{087})--(\ref{088}) with $\bmu(\bX \hbbeta_n) = (e^{\theta_1}/(1 + e^{\theta_1}),
\cdots, e^{\theta_n}/(1 + e^{\theta_n}))\t$, $\Sig(\bX \hbbeta_n) = \diag\{e^{\theta_1}/(1 + e^{\theta_1})^2, \cdots, e^{\theta_n}/(1 + e^{\theta_n})^2\}$, and $(\theta_1, \cdots, \theta_n)\t = \bX \hbbeta_n$.

\subsection{Poisson regression} \label{secB.3}
In the Poisson regression model, $b(\theta) = e^\theta$, $\theta \in \mathbb{R}$ and the quasi-log-likelihood of the sample is $\ell_n(\by, \bbeta) = \by\t \bX \bbeta - \bone\t \bb(\bX \bbeta)$ with the last constant term in (\ref{002}) ignored, which attains its maximum at $\hbbeta_n$. Since
$b'(\theta) = e^\theta$ and $b''(\theta) = e^\theta$, matrices $\widehat{\bA}_n$ and $\widehat{\bB}_n$ are given in (\ref{087})--(\ref{088}) with
$\bmu(\bX \hbbeta_n) = (e^{\theta_1}, \cdots,
e^{\theta_n})\t$ and $\Sig(\bX \hbbeta_n) =
\diag\{e^{\theta_1}, \cdots, e^{\theta_n}\}$, and $(\theta_1, \cdots, \theta_n)\t = \bX \hbbeta_n$.

\end{document}